\documentclass{siamart220329}

\usepackage[numbers,sort&compress]{natbib}  

\usepackage{graphicx}

\usepackage{listings}
\definecolor{mygreen}{RGB}{28,172,0} 
\definecolor{mylilas}{RGB}{170,55,241}
\usepackage[utf8]{inputenc}

\usepackage{flexisym}	 
\usepackage{breqn}       

\usepackage{amsmath,amssymb}
\usepackage{mathtools}   
\usepackage{physics}
\usepackage{bm} 

\usepackage{url}

\usepackage[utf8]{inputenc} 

\usepackage{csquotes}
\MakeOuterQuote{"}       

\usepackage[colorinlistoftodos]{todonotes}   

\usepackage{algorithm,algpseudocode}

\usepackage{dirtree}
\usepackage{framed}

\usepackage{subcaption}
\numberwithin{equation}{section}
\numberwithin{figure}{section}

\usepackage{xr}

\DeclarePairedDelimiterX{\innerp}[2]{\langle}{\rangle}{#1,#2}
\renewcommand{\ip}[2]{\innerp*{#1}{#2}}

\newcommand{\normp}[2]{{\norm{#1}}_{#2}}

\newcommand{\Tee}{\mathsf{T}} 
\newcommand{\transpose}[1]{{#1}^\Tee}

\newcommand{\FF}{\mathcal{F}}
\newcommand{\VV}{\mathcal{V}}
\newcommand{\EE}{\mathcal{E}}
\newcommand{\vertex}{\mathtt{v}}
\newcommand{\edge}{\mathtt{e}}
\newcommand{\Laplacian}{\triangle}
\newcommand{\degree}{d}
\DeclareMathOperator{\diag}{diag}

\newcommand{\I}{\vb{I}}

\renewcommand{\P}{\vb{P}}
\newcommand{\A}{\vb{A}}
\newcommand{\Lplus}{\vb{L}_+}
\newcommand{\Lminus}{\vb{L}_-}

\newcommand{\Lint}{\vb{L}_{\mathrm{int}}}
\newcommand{\LVC}{\vb{L}_{\mathrm{VC}}}
\newcommand{\LZero}{\vb{L}_0}

\newcommand{\PVC}{\vb{P}_{\mathrm{VC}}}
\newcommand{\PZero}{\vb{P}_0}
\newcommand{\Pint}{\vb{P}_{\mathrm{int}}}

\newcommand{\MVC}{\vb{M}_{\mathrm{VC}}}
\newcommand{\MBC}{\vb{M}_{\mathrm{BC}}}
\newcommand{\MNH}{\vb{M}_{\mathrm{NH}}}

\newcommand{\M}{\vb{M}}
\newcommand{\Nint}{N_{\mathrm{int}}}
\newcommand{\Next}{N_{\mathrm{ext}}}
\newcommand{\Fint}{\mathbb{F}^{\mathrm{int}}}
\newcommand{\Fext}{\mathbb{F}^{\mathrm{ext}}}
\newcommand{\DD}{\vb{D}^{2}}

\renewcommand{\SS}{\vb{S}}

\renewcommand{\u}{\vb{u}}

\newcommand{\x}{\vb{x}}
\newcommand{\xInt}{\x^{\mathrm{int}}}
\newcommand{\xExt}{\x^{\mathrm{ext}}}
\newcommand{\xkInt}[1]{x_{#1}^{\mathrm{int}}}
\newcommand{\xkExt}[1]{x_{#1}^{\mathrm{ext}}}

\newcommand{\f}{\vb{f}}

\newcommand{\e}{\vb{e}}

\renewcommand{\P}{\vb{P}}

\newcommand{\D}{\vb{D}}
\newcommand{\Zero}{\vb{0}}

\newcommand{\vphi}{\boldsymbol{\phi}}
\newcommand{\vpsi}{\boldsymbol{\psi}}

\newcommand{\stackmat}[2]{
\left(   \begin{array}{c}
#1 \\ \hline #2
\end{array} \right) 
}

\newcommand{\RR}{{\mathbb R}}

\DeclarePairedDelimiter{\set}{\{}{\}}

\newcommand{\Eye}{\vb{I}}
\newcommand{\matstack}[2]{\left( \begin{array}{c} #1 \\ \hline  #2 \end{array} \right)}

\newcommand{\qGraph}{\lstinline|quantumGraph|} 
\newcommand{\qg}{\lstinline|qg|}

\newcommand{\eigs}{\lstinline|eigs|}
\newcommand{\digraph}{\lstinline|digraph|}
\newcommand{\ttG}{\lstinline|G|}
\newcommand{\source}{\lstinline|source|}
\newcommand{\target}{\lstinline|target|}
\newcommand{\Uniform}{\lstinline|'Uniform'|}
\newcommand{\Chebyshev}{\lstinline|'Chebyshev'|}
\newcommand{\none}{\lstinline|'None'|}
\newcommand{\NaN}{\lstinline|NaN|}
\newcommand{\ttL}{\lstinline|L|}
\newcommand{\tty}{\lstinline|y|}
\newcommand{\ttx}{\lstinline|x|}
\newcommand{\ttnx}{\lstinline|nx|}
\newcommand{\Edge}{\lstinline|Edge|}

\newcommand{\listtt}[1]{\lstinline|#1|}

\renewcommand\descriptionlabel[1]{$\bullet$ \texttt{#1}}
\newcommand{\itemtt}[1]{\item{\listtt{#1}}}

\renewenvironment{framed}[1][\hsize]
  {\MakeFramed{\hsize#1\advance\hsize-\width \FrameRestore}}%
  {\endMakeFramed}

\makeatletter

\newcommand*{\addFileDependency}[1]{
\typeout{(#1)}
\@addtofilelist{#1}
\IfFileExists{#1}{}{\typeout{No file #1.}}
}\makeatother

\ifpdf
\hypersetup{ pdftitle={QGLAB: A MATLAB Package for Computations on Quantum Graphs} }
\fi
\newcommand{\suppendix}{appendix}

\graphicspath{{Images/}}

\title{QGLAB: A MATLAB Package for Computations on Quantum Graphs\thanks{\funding{R.H.G. acknowledges support from the NSF through NSF Grant DMS-2206016 G.C. and J.L.M.
acknowledge support from the NSF through NSF CAREER Grant
DMS-1352353, NSF Grant DMS-1909035 and NSF FRG grant DMS-2152289.}}}

\author{Roy H. Goodman\thanks{Department of Mathematical Sciences, New Jersey Institute of Technology, Newark, NJ, (\email{goodman@njit.edu}).}
\and Grace Conte\thanks{Johns Hopkins University Applied Physics Laboratory, Baltimore, MD, (\email{gracie.conte@jhuapl.edu}).}
\and Jeremy L. Marzuola\thanks{Department of Mathematics, University of North Carolina, Chapel Hill, NC (\email{marzuola@email.unc.edu}).}
}
\headers{QGLAB: A MATLAB Package for Computations on Quantum Graphs}{Roy H. Goodman, Grace Conte, and Jeremy Marzuola}

\date{\today}            

\begin{document}
\maketitle

\begin{abstract}
We describe QGLAB, a new MATLAB package for analyzing partial differential equations on quantum graphs.  The software is built on the existing, object-oriented MATLAB directed-graph class, inheriting its structure and adding additional easy-to-use features.  The package allows one to construct a quantum graph and accurately compute the spectrum of elliptic operators, solutions to Poisson problems, the linear and nonlinear time evolution of a variety of PDEs, the continuation of branches of steady states (including locating and switching branches at bifurcations) and more. It overcomes the major challenge of discretizing quantum graphs---the enforcement of vertex conditions---using non-square differentiation matrices. It uses a unified framework to implement finite-difference and Chebyshev discretizations of differential operators on a quantum graph. For simplicity, the package overloads many built-in MATLAB functions to work on the class.
\end{abstract}

\section{Introduction}

This paper introduces the main ideas used to build QGLAB, a software package written in MATLAB for computations on quantum graphs, and provides several examples of its use and accuracy~\cite{Goodman:2024}. The supplementary materials present more thorough operating instructions and additional computational examples. 

Quantum graphs, networks of one-dimensional edges interacting via vertex conditions, appear in the literature going back many years. The modern study of the subject begins with an analysis of their spectral statistics in Ref.~\cite{kottos1997quantum}, whose authors coined the term "quantum graphs."  The spectral theory and properties of quantum graph operators were further developed in Ref.~\cite{exner2005convergence}, where quantum graphs were realized as the limits of quantum equations on thin wire-like domains; see also~\cite{grieser2008spectra,exner2013general} and the references therein. Quantum graphs provide effectively one-dimensional model equations that enable explicit calculations that serve as a backbone for representing geometric and spectral theoretic properties of more complicated higher-dimensional quantum models.  For further introduction and history, we recommend Refs.~\cite{berkolaiko2017elementary, Berkolaiko:2013bq}.

Numerical packages are essential to facilitate further study for several common reasons: making progress on larger-scale problems and those with nonlinearities and time dependence all depend on intuition built from numerical experimentation and visualization. This project provides high-level tools that allow users to quickly and easily set up, solve, and visualize the solutions to problems posed on quantum graphs. It overloads many built-in MATLAB commands for basic calculations and plotting.

The package's foundation is a quantum graph class built on MATLAB's directed-graph class. While we have striven to write a general-purpose software package for quantum graph computations, the direction of development has been guided by two classes of problems of research interest to its authors:
\begin{enumerate}
\item Computing bifurcation diagrams: standing waves of NLS occur along one-parameter families or \emph{branches} rather than at isolated points. Such branches must be computed using continuation methods to understand the solutions' parameter dependence. Branches may cross, and the stability of solutions change at isolated \emph{bifurcation points}. The package can detect the most common bifurcations and switch branches; see Refs.~\cite{Berkolaiko:2020, Goodman:2019go}.

\item Spectral accuracy in space and high-order time stepping: a planned future project is to compute time-periodic and time-relative-periodic orbits of the full time-dependent NLS on a compact quantum graph.
\end{enumerate}

QGLAB grew from numerical studies in the authors' research~\cite{Goodman:2019go, Kairzhan:2019fe, Marzuola:2016bl}. Other works have developed numerical packages for quantum graphs. The most complete is \emph{GraFiDi}, a Python-based package described in~\cite{Besse:2021}, which overlaps with QGLAB but has fewer features. Malenova built a small quantum graph package using \emph{Chebfun}~\cite{malenova2013spectra}. Others have studied finite-difference, finite-element, and Galerkin-based numerical methods without creating a software package for general use~\cite{Arioli:2017dh, Berkolaiko:2020, Bottcher:2024, brio2022spectral}. 

QGLAB is distinguished by the breadth of problems it can handle: linear and nonlinear, stationary and time-dependent, constant and variable coefficient. It separates the solution of differential equations along a graph's edges from the constraint of solving the vertex conditions. This has allowed us to develop a language in which the numerical algorithms for solving problems on quantum graphs can be expressed in just a few lines, as demonstrated in numerous examples.

\subsection{Defining a quantum graph}
A quantum graph \(\Gamma\) consists of a directed metric graph, considered as a \emph{complex of edges}, on which a function space and differential operators are defined. To be more specific, we define the graph \(\Gamma = \left(\VV,\EE \right) \) as a set of vertices \( \VV = \set*{\vertex_n, \, n= 1,\ldots,\abs*{\VV}} \) and a set of directed edges \( \EE=\set*{\edge_m=(\vertex_i \to \vertex_j),\, m= 1,\ldots,\abs*{\EE}} \) and to each edge assign a positive length \(\ell_m\) and impose upon the edge a coordinate \(x\) that increases from \(0\) to \(\ell_m\) as the edge is traversed from the \(\vertex_i\) to \(\vertex_j\). In general, we may consider the lengths of some or all of the edges infinite, in which case that edge will be connected to a single vertex, but we defer this to later publications.  Define the \emph{degree} \(\degree_n\)of vertex \(\vertex_n\) as the number of edges that include that vertex as an initial or final point, counting twice if an edge connects the vertex to itself.

\begin{figure}[htbp] 
   \centering
   \includegraphics[width=2in]{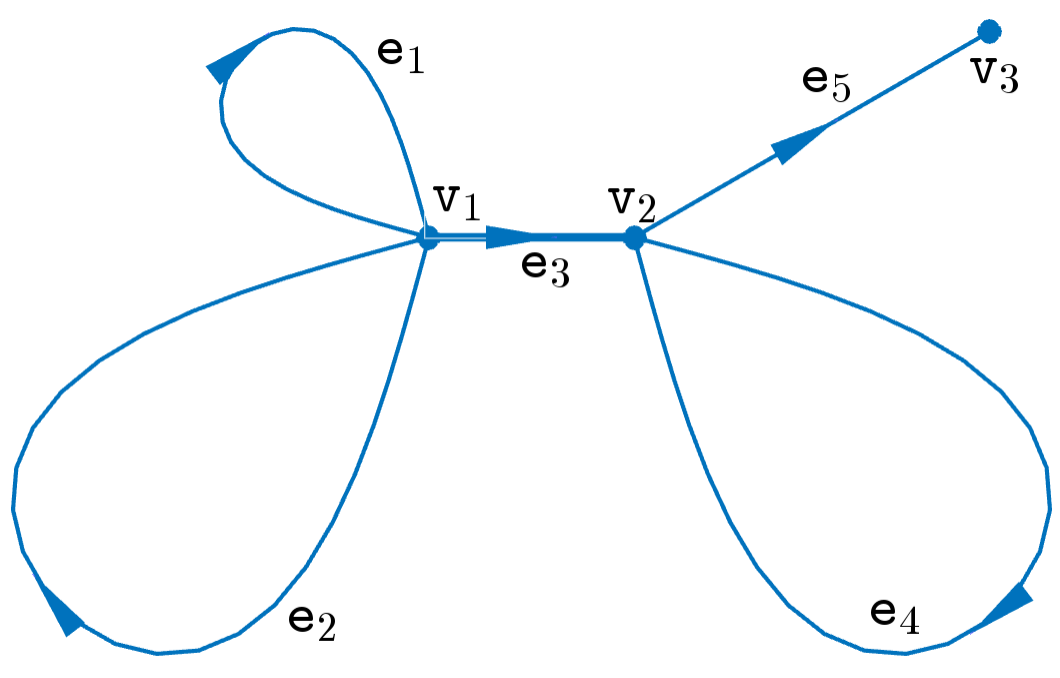} 
   \caption{A directed graph with three vertices and five edges.}
   \label{fig:sillyGraph}
\end{figure}

An graph with \(\abs*{\VV}=3\) vertices and \(\abs*{\EE}=5\) edges is shown in Fig.~\ref{fig:sillyGraph}. The  vertices have degrees \(\degree_1=5\), \(\degree_2=3\), and \(\degree_3=4\). Under our definition, multiple edges may share the same initial and final vertex, as do edges \(\edge_1\) and \(\edge_2\).

A function \(\Psi(x)\) defined on \(\Gamma\) is defined as a collection of functions on each of the edges \(\left.\Psi\right|_{\edge_m} = \psi_m(x)\). The Laplace operator on the graph is defined by%
\footnote{With the inclusion of such a potential, the operator $-\Laplacian$ is more properly called a Schrödinger operator, but we abuse terminology for simplicity.} 
\begin{equation} \label{Laplace_V}
\left. \Laplacian \right\rvert_{\edge_m} =  \dv[2]{x} -V_m(x), \qfor {0<x<\ell_m},
\end{equation}
subject to appropriate compatibility conditions at the vertices, given the presence of a potential $V(x)$ with $\left. V(x)\right\rvert_{\edge_m}=V_m(x)$ which we usually set to zero. To define a Laplacian requires a function space. The graph and vertex conditions define a quantum graph, and we are most concerned with vertex conditions giving rise to a self-adjoint operator. The text~\cite{Berkolaiko:2013bq} describes general criteria for self-adjoint vertex conditions. We discuss some of these here.

We assume the function is continuous at each vertex, yielding \(\degree_n -1\) conditions at vertex  \( \vertex_n \). Let \( \VV_n \) be the set of edges adjacent to this vertex, double counting self-directed edges. Continuity defines a function value at the vertex to be
\begin{equation}
\Psi(\vertex_n) \equiv \psi_i (\vertex_n) = \psi_j (\vertex_n), \forall \edge_i, \edge_j \in \VV_n.
\label{continuity}
\end{equation}
We define the weighted Robin-Kirchhoff vertex condition as
\begin{equation}
\sum_{\edge_m \in \VV_n} w_m \psi_m'(\vertex_n) + \alpha_n \Psi(\vertex_n) = 0, 
\label{KirchhoffRobin}
\end{equation}
where the derivative is, in all cases, taken in the direction \emph{pointing away} from the vertex. In the case \(w_m \equiv 1\), \(\alpha_n\equiv0\), this reduces to the Neumann-Kirchhoff vertex condition, which is the natural generalization of the Neumann boundary condition on a line segment. Interpreting the equation \(\Laplacian \Psi=0\) as describing a steady state of the heat equation, the Neumann-Kirchhoff vertex condition states that the net heat flux into the vertex vanishes. Associating non-unit weights \(w_m\) to the edges then specifies that the flux from each edge is proportional to its weight. Such weights define \emph{balanced} star graphs in Ref.~\cite{Kairzhan:2019fe}. Setting \(\alpha_n \neq 0 \)  generalizes a Robin vertex condition and may be interpreted as a delta-function potential at the vertex. The Dirichlet condition is 
\begin{equation}
\Psi(\vertex_n) = 0.
\label{Dirichlet}
\end{equation}
Setting the right-hand side of~\eqref{KirchhoffRobin} or~\eqref{Dirichlet} to a value \( \phi_n \neq 0\) defines
nonhomogeneous vertex conditions.

We define norms and function spaces, e.g., \(L^p(\Gamma)\):
\begin{equation} \label{Lpnorm}
\normp{\Psi}{L^p(\Gamma)}^p = \sum_{m=1}^{\abs*{\EE}} w_m \normp{\psi_m}{L^p}^p
\end{equation}
and the \(L^2\) inner product
\begin{equation}
\ip{\Psi}{\Phi} = \sum_{m=1}^{\abs*{\EE}} w_m \int_0^{\ell_m} \psi_m^*(x) \phi_m(x) \dd{x}.
\label{innerProduct}
\end{equation}
The weights \(w_m\) must appear in these definitions for the conservation laws for evolution equations below and to make the Laplace operator~\eqref{Laplace_V} self-adjoint. 
We define \(H^1(\Gamma)\) as the space of square-integrable functions with square-integrable first derivatives. More subtly, we define \(L^2(\Gamma)\), \(H^1(\Gamma)\), and \(H^2(\Gamma)\) to be the space of functions that are in each of these function spaces edgewise, but we define \(L^2_\Gamma\) to be the space \(L^2(\Gamma)\) equipped with the inner product~\eqref{innerProduct}. Similarly, \(H^1_\Gamma\) is the space of functions in \(H_1(\Gamma)\) satisfying the continuity condition~\eqref{continuity}, and \(H^2_\Gamma\) consists of functions in \(H^2(\Gamma)\) satisfying both Eq.~\eqref{continuity} and either the vertex condition~\eqref{KirchhoffRobin} or~\eqref{Dirichlet}.

\subsection{The Eigenvalue Problem}
\label{sec:eigenvalueProblem}
The first natural question about the Laplacian operator defined on \(\Gamma\) is its spectrum and eigenfunctions. Such properties have  been studied  extensively; for a small subset of relevant works, see for instance~\cite{alon2018nodal,band2012number,Berkolaiko:2013bq,berkolaiko2017elementary,berkolaiko2019surgery,harrell2016spectral} and the references therein. 

The compact quantum graph Laplacian has only discrete spectrum, so we must compute the set of eigenvalues \(\lambda\) and eigenfunctions \(\Psi\) such that 
\begin{equation}
    \Laplacian \Psi = \lambda \Psi.
    \label{analyticEvalProb}
\end{equation}
The spectrum is countable, unbounded below, and has finitely many positive values. The non-positive eigenvalues \(\lambda=-k^2\) can be found by seeking the analytic solution
\begin{equation}\label{edgesoln}
    \psi_m(x) = a_me^{ikx} + b_me^{ik(\ell_m-x)} \qquad m=1,2,\ldots,|\EE|.
\end{equation}
The vertex conditions form a homogeneous system of \(2\abs{\EE}\) linear equations. Its solution requires the vanishing of the determinant of the associated matrix \(\SS(k)\), a function \(\Sigma(k)\) called the \emph{secular determinant} which can be normalized to take real values when \(k\in \RR\)~\cite{Berkolaiko:2013bq}.
The recent dissertation~\cite{Conte:2022} shows this holds under the more general vertex conditions~\eqref{KirchhoffRobin}. In addition to QGLAB's many numerical features, it can symbolically compute the graph's real-valued secular determinant.

Two recent publications note numerically computing a determinant and finding its zeros have high complexity and low accuracy, so that finding the zeros of $\Sigma$ is a poor method for computing the spectrum~\cite {Bottcher:2024,brio2022spectral}. Both suggest that finding the values $k^*$ where the condition number of $\SS(k)$ diverges is faster and more accurate; they subsequently find the eigenfunction as null vectors of \(\SS(k^*)\). Further, both suggest using the null vectors of $\SS(k^*)$ as the basis for Galerkin methods to solve stationary and evolutionary PDE on metric graphs. This approach is ill-suited to nonlinear problems and those with edgewise-defined potentials or semi-infinite edges: nonlinearities are cumbersome to implement in Galerkin methods, eigenfunctions in the presence of potentials do not take the form~\eqref{edgesoln}, and graphs with infinite edges may not even have point spectrum. The present work demonstrates numerical solutions to the first two problems. Extension to infinite edges is planned.

\subsection{PDE on a Quantum Graph}

Our primary motivating problem for building QGLAB is the nonlinear Schrödinger equation
\begin{equation}
i \pdv{\Psi}{t} = \Laplacian \Psi + (\sigma+1) \abs*{\Psi}^{2\sigma} \Psi
\label{NLS},
\end{equation}
where \(\sigma \ge 0\) and \(\sigma=1\) is the most commonly-studied cubic case.
We are especially interested in the stationary NLS obtained by assuming \(\Psi(x,t) = e^{i \Lambda t} \Psi(x) \),
\begin{equation}
\FF(\Psi,\Lambda) \equiv 
\Lambda \Psi + \Laplacian \Psi + (\sigma+1) \abs*{\Psi}^{2\sigma} \Psi = 0.
\label{stationary_NLS}
\end{equation}
The evolution of Eq.~\eqref{NLS} conserves both the \(L^2\) norm defined in Eq.~\eqref{Lpnorm} and an energy
\begin{equation} \label{energyNLS}
    E(\Psi) =    \normp{\Psi'}{L^2(\Gamma)}^2 
                - \normp{\Psi}{L^{2(\sigma+1)}_\Gamma}^{2(\sigma+1)}
                + \sum_{\edge_m \in \EE} w_m\int_{\edge_m} V_m(x) \abs{\Psi_m(x)}^2 \dd x
                + \sum_{\vertex_n \in \VV} \alpha_n \abs*{\Psi(\vertex_n)}^2 ,
\end{equation}
where \(\Psi'\) is defined edge-by-edge. The NLS equation on the real line also conserves a momentum functional, but NLS on quantum graphs does not unless certain other restrictions to the weights and initial conditions hold; see~\cite{Kairzhan:2019fe}. 

Linear and nonlinear PDEs on quantum graphs are longstanding and active subjects. Many groups have studied Eq.~\eqref{stationary_NLS} as surveyed in Refs.~\cite{noja2014nonlinear,borrelli2019overview}. The existence of ground states and the solution stability of stationary solutions are two important questions often studied from a variational perspective.  There are too many works in this direction to properly do the subject justice, but see for instance~\cite{adami2017nonlinear, cacciapuoti2018variational, de2023notion,  dovetta2020ground, borrelli2019nonlinear, serra2016lack} and the references within to get a sense of the field.  Others have studied the existence and stability of stationary states for Dirac and KdV equations~\cite{borrelli2019nonlinear,pava2021linear}.  Recent works, including~\cite{Berkolaiko:2020,noja2015bifurcations,gnutzmann2016stationary, Marzuola:2016bl, Goodman:2019go}, use asymptotic and bifurcation-theoretical approaches to analyze the existence of multiple branches of solutions to~\eqref{stationary_NLS}.  The book~\cite{mugnolo2014semigroup} gives an excellent introduction to time-dependent PDE on graphs. References including~\cite{Kairzhan:2019fe, mugnolo2018airy, sabirov2018dynamics} analyze the time-dependent phenomena exhibited by  Schrödinger, Dirac, and KdV equations on graphs.  This short overview gives a flavor of the questions that can be posed on quantum graphs and the breadth of topics yet to be explored.

QGLAB has been explicitly written for PDE with Laplacian spatial derivative terms. In addition to the previously discussed NLS equation, these include the wave equation (\cite{noja2014nonlinear}), heat equation (\cite{becker2021schrodinger,borthwick2022heat}), and their nonlinear cousins such as the nonlinear Klein-Gordon equations, including sine-Gordon~\cite{goloshchapova2021nonlinear, maier2022breather, sabirov2018dynamics, sobirov2016sine} and the  Kolmogorov–Petrovsky–Piskunov (KPP) equation~\cite{du2020fisher}, all of which can be defined on a quantum graph. QGLAB provides examples of solving all of these PDEs.

\subsection{Organization of the paper}

Section~\ref{sec:NumericalMethods} discusses the numerical methods QGLAB uses to discretize and solve equations posed on quantum graphs. 
The longest part, Subsection~\ref{sec:discretization}, discusses the overall framework of the discretization and its implementation using both finite-difference and Chebyshev approximations of derivatives and the implementation of the vertex conditions.  
We apply this framework to discretize eigenvalue problems in Subsection~\ref{sec:numerical_eigenproblems}, where we also discuss the symbolic calculation of the secular determinant for the general class of vertex conditions discussed. 
Subsection~\ref{sec:solvers} describes the nonlinear solvers and continuation algorithms, while Subsection~\ref{sec:timestepping} describes the implementation of time steppers for evolution equations. 
Section~\ref{sec:MATLAB} is devoted to the MATLAB implementation of the tools discussed in QGLAB, including a discussion of MATLAB's directed-graph class in Section~\ref{sec:digraph} and the QGLAB's quantum graph class, which is built on top of this, in Section~\ref{sec:initializing}. 
Section~\ref{sec:operations} discusses basic operations on class objects.
In Section~\ref{sec:examples}, we consider three examples illustrating QGLAB in practice: a Poisson problem, an eigenvalue problem, and an initial-value problem for the cubic NLS equation. 
We summarize our contributions and give an outlook on potential future features and applications in Section~\ref{sec:conc}.  
The \suppendix\ contains extensive additional materials in two sections. The first, Sec.~\ref{sec:appendix_examples}, is devoted to demonstrating both the implementation and efficacy of QGLAB on a variety of examples, including stationary problems---eigenvalue problems, the Poisson equation, and the computation and continuation of standing waves and evolutionary PDE problems. The second, Sec.~\ref{sec:function_listing}, contains a complete listing of user-callable function definitions and explicit instructions for their use.

\section{Numerical Methods}
\label{sec:NumericalMethods}

This section discusses the numerical methods used to implement the quantum graph class and solve various problems. It briefly presents some examples described in detail in Sec.~\ref{sec:examples} of the \suppendix.

\subsection{Discretization and vertex conditions}
\label{sec:discretization}

QGLAB can perform tasks including solving for nonlinear standing waves and numerically integrating evolution equations on a quantum graph. Most important is discretizing the Laplace operator and solving the Laplace and Poisson equations. It provides two discretizations: centered differences and Chebyshev collocation. They can be used interchangeably, although the Chebyshev discretization is, by construction, more accurate. 

The two discretizations are implemented using a common framework:  the function $\Psi(x)$ is approximated on an \emph{extended grid} \(\xExt\) with enough points to approximate both the PDE solution and the vertex conditions. In contrast, the discretized PDE is satisfied on a smaller \emph{interior} grid \(\xInt\) containing two fewer points per edge.
Thus, the discrete Laplacian matrix is non-square, with $2\abs{\EE}$ more columns than rows, mapping from approximations on \(\xExt\) to approximations on \(\xInt\).  The vertex conditions are implemented as constraints, not incorporated directly into the discretized Laplacian matrix. This choice has several attractive features described below.

Driscoll and Hale introduced non-square differentiation matrices using rectangular discretization matrices for use in the Chebfun package~\cite{driscoll2016rectangular,chebfun:2014}. Aurentz and Trefethen have written an excellent review, developing the theory for the \emph{block operators} that implement these ideas via a sequence of well-chosen examples~\cite{Aurentz:2017}.

\subsubsection{Finite-difference discretization}
\label{sec:finitediff}

\begin{figure}[htbp] 
   \centering
   \includegraphics[width=5in]{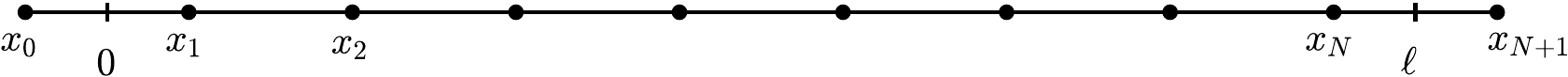}
   \caption{Discretization of the interval \([0,\ell]\) using ghost points.}
   \label{fig:ghostPointSchematic}
\end{figure}

The finite-difference method is implemented using second-order centered differences with the boundary conditions enforced at so-called ghost points, as discussed, for example, in the textbook~\cite[Sec.\ 4.2.2]{Edsberg:2016wu}. We review this technique for discretizing the two-point boundary value problem
\begin{equation}
\label{robinSegment}
\dv[2]{u}{x} -V(x) u(x)= f(x), 0<x<\ell, \; 
u'(0)  + \alpha_0 u(0) = \phi_0, \; 
-u'(\ell)  + \alpha_\ell u(\ell) = \phi_\ell 
\end{equation} 
and then discuss the straightforward extension to quantum graphs. The sign on the \(u'\) term in the boundary conditions is chosen to agree with our quantum graph convention that all derivatives are taken in the direction pointing away from a vertex of the quantum graph in the definition of vertex conditions. Given a discretization length \(h = \frac{\ell}{N}\), place points on a grid offset by half a step size \(x_k = \left(k-\frac{1}{2} \right) h\) for \(0\le k \le N+1\) as shown in Fig.~\ref{fig:ghostPointSchematic}. The first and last points lie outside the interval of interest, and the endpoints of the desired interval do not appear on the list of points. Letting \(u_k\) approximate \(u(x_k)\), the discretized equation at the interior points is then
\[
u''(x_k) - V(x_k)u(x_k) \approx \frac{u_{k-1} - 2 u_k + u_{k+1}}{h^2}  =f_k = f(x_k), \qfor k = 1,\ldots, N,
\]
up to an error of \(\order{h^2}\). The value of $u(0)$ and $u'(0)$ may be approximated to second order via suitable linear combinations of $u(\pm h/2)$, yielding a second order approximation of the boundary condition Eq.~\eqref{robinSegment}:
\begin{equation}
\left(\frac{\alpha_0}{2} -\frac{1}{h} \right) u_0 +
\left(\frac{\alpha_0}{2} +\frac{1}{h} \right) u_1 = \phi_0
\label{discreteBC}
\end{equation}
and a similar approximation for the right boundary condition. In the case of a Dirichlet boundary condition at \(x=0\), this is replaced by \( \frac{1}{2} \left(u_0 +  u_1\right) = \phi_0\). 

Assembling the equations for the second derivatives and the boundary conditions into matrix--vector form, we define the vectors 
\begin{equation}
\u = \transpose{\left(u_0,u_1,\ldots,u_{N+1} \right)}, 
\f = \transpose{\left(f_0,f_1, \ldots, f_{N+1} \right)}, 
\qand
\vphi = \transpose{\left( \phi_0, \phi_\ell \right)},
\label{u_f_phi}
\end{equation}
as well as the \(N\times(N+2)\) \emph{interior projection matrix}
\begin{equation}
\Pint = \begin{pmatrix}
\Zero_{N \times 1} & \Eye_N & \Zero_{N \times 1} 
\end{pmatrix},
\label{Pint}
\end{equation}
the \(N\times(N+2)\) \emph{interior Laplacian matrix}
\begin{equation}
\Lint = \frac{1}{h^{2}}
\begin{pmatrix}
1 & -2 & 1 &  & & \\
 & 1 & -2 & 1 &  &\\
 &  & \ddots  & \ddots &\ddots  & & \\
  &  &  & 1 &-2 & 1
\end{pmatrix}
- \Pint \cdot \diag{V(x_k)}
,
\label{Lint}
\end{equation}
and the \(2\times(N+2) \) \emph{boundary condition matrix}
\begin{equation}
\MBC = 
\begin{pmatrix}
\left( \tfrac{\alpha_{0}}{2}-\tfrac{1}{h}\right)  & \left( \tfrac{\alpha_{0}}{2}+\tfrac{1}{h}\right)  & \cdots & 0 & 0 \\
0 & 0 & \cdots & \left( \tfrac{\alpha_{\ell}}{2}+\tfrac{1}{h}\right)  & \left( \tfrac{\alpha_{\ell}}{2}-\tfrac{1}{h}\right) 
\end{pmatrix},
\label{MBC}
\end{equation}
where \(\Eye_N\) is the \(N\)-dimensional identity matrix and \(\Zero_{M\times N}\) is a matrix of size \((M\times N)\) of all zeros. The matrices \(\Lint\) and \(\Pint\) are linear maps from  approximation of $u$  on the \emph{extended grid} \(\xExt =\left\{ x_0,\ldots,x_{N+1} \right\}\) to its approximation on the \emph{interior grid} \(\xInt = \left\{ x_1,\ldots,x_{N} \right\}\). 
With this, we discretize the differential equation as
\begin{equation}\label{discSecondDeriv}
    \Lint \u = \Pint \f,
\end{equation}
and the boundary conditions as 
\begin{equation}\label{discBCs}
    \MBC \u = \vphi.
\end{equation}
We can then combine these into a single system of \( N+2 \) equations in \(N+2\) unknowns
\begin{equation}\label{discSystemEqs}
\matstack{\Lint}{\MBC} \u = 
\matstack{\Pint}{\Zero_{2\times(N+2)}} \f + 
\matstack{\Zero_{N\times 2}}{\Eye_{2}} \vphi.
\end{equation}

We make two brief remarks on this approach. First, it is more common to solve the discretized boundary condition~\eqref{discBCs} for \(u_0\) and \(u_{N+1}\)~\cite{brio2022spectral}. Inserting these values into Eq.~\eqref{discSecondDeriv}  yields a reduced system with \(N\) unknowns. We leave the system in form~\eqref{discSystemEqs} for two reasons: first, it makes implementing non-homogeneous boundary conditions very easy, and the other is that it makes the approach more similar to how we handle boundary conditions using Chebyshev discretization. Second, the null space of the matrix \(\Lint\) mimics that of the second derivative: it consists of vectors \(\vb{v}\) with \(v_n = a n + b\) and is two-dimensional. 

This scheme extends easily to the Poisson problem on the quantum graph,
\begin{subequations}
\begin{gather}
\Laplacian \Psi(x)  = f(x), \qqtext{i.e.,} \psi_m''(x) -V_m(x) \psi_m(x)= f_m(x), \edge_m, 1 \le m \le \abs*{\EE}\\
\psi_i (\vertex_n) = \psi_j (\vertex_n), \forall \edge_i, \edge_j \in \VV_n, \qqtext{i.e., continuity,} \label{continuitycondition}\\
\sum_{\edge_m \in \VV_n} w_m \psi_m'(\vertex_n) + \alpha_n \Psi(\vertex_n) = \phi_n 
\qor
\Psi(\vertex_n) = \phi_n, 1 \le n \le \abs*{\VV}. \label{VC}
\end{gather}
\label{Poisson}
\end{subequations}
We discretize each edge \(\edge_m\) with the ghost-point formulation, with \(N_m + 2\) discretization points  defining \(\xExt_m\), and a mesh size \(h_m = \ell_m/N_m\), generating \((N_m)\times N_m+2\) matrices \(\Lint^{(m)}\) and \(\Pint^{(m)}\)of the same forms as matrices~\eqref{Lint} and~\eqref{Pint}. Thus, letting \(\Nint= \sum_{m=1}^{\abs*{\EE}} N_m\)  and \(\Next=\Nint+2\abs*{\EE}\), this results in \(\Next\) unknowns arranged as 
\begin{equation}\label{vpsi}
\vpsi = \begin{pmatrix} \vpsi^{(1)} \\ \vdots \\ \vpsi^{(\abs*{\EE})} \end{pmatrix},
\qqtext{where}
\vpsi^{(m)} = \begin{pmatrix} \vpsi^{(m)}_0 \\ \vdots  \\ \vpsi^{(m)}_{N_m+1} \end{pmatrix}.
\end{equation}
The vector \(\f\) is assigned similarly, and the vector of nonhomogeneous boundary terms is \( \vphi = \transpose{\left(\phi_0,\ldots,\phi_{\abs*{\VV}} \right)} \).
 Enforcing the continuity condition~\eqref{continuitycondition} at the vertex \(\vertex_n\) requires \((d_n-1)\) rows and the Robin-Kirchhoff condition~\eqref{VC} at vertex \(\vertex_n\) involves the \(2 d_n\) adjacent discretization points in one row. The derivative and function values at the vertex are approximated to second order using a straightforward generalization of the reasoning leading to Eqs.~\eqref{discreteBC} and~\eqref{MBC}.
Altogether, these form a matrix \(\MVC^{(n)}\) of dimension \((2 d_n) \times \Next\).
We let 
\begin{equation} \label{Luniform}
\LVC = 
\stackmat{\Lint}{\MVC} 
=
\left(
\begin{array}{c}
\mqty{\dmat{\Lint^{(1)},\ddots,\Lint^{\left(\abs*{\EE}\right)}}} \\
\hline
  \MVC^{(1)}  \\
  \vdots  \\
  \MVC^{\left(\abs*{\VV}\right)}
\end{array}
\right), 
\end{equation}
and
\begin{equation} \label{Puniform}
\PZero = 
\stackmat{\Pint}{\Zero_{2\abs*{\EE} \times \Next}}
=
\left(
\begin{array}{c}
\mqty{\dmat{\Pint^{(1)},\ddots,\Pint^{\left(\abs*{\EE}\right)}}} \\
\hline
\Zero_{2\abs*{\EE} \times (\Next)}
\end{array}
\right).
\end{equation}
We define the \emph{nonhomogeneity matrix} \(\MNH\) in two steps. First define a matrix \(\M\) of size  \( 2\abs{\EE} \times \abs{\VV} \)
such that 
\[
\M(i,j) = \begin{cases}
1 & \text{If the \(j\)th Neumann-Kirchhoff condition is enforced by row \(i\) of \(\M\)}, \\
0 & \text{otherwise}
\end{cases}
\]
and then define
\begin{equation} \label{MNH}
\MNH = \left(
\begin{array}{c}
\Zero_{\Nint,\abs{\VV}} \\
\hline
\M
\end{array}
\right).
\end{equation}
This assigns each entry of the nonhomogeneous term \(\vphi\) to the row of the system enforcing the vertex conditions. The discretization of system~\eqref{Poisson} can now be represented as
\begin{equation}
\LVC \vpsi   = \PZero \f + \MNH \vphi.
\label{discretizedLaplacian} \\
\end{equation}
The use of this discretization is demonstrated in Sec.~\ref{sec:poisson_example}.

We introduce some notation to simplify our discussion of the numerical problems addressed below. The matrices \(\Pint\) and \(\Lint\), of dimension \(\Nint \times \Next\), represent linear maps from the function space \(\Fext\), of functions defined on the extended grid \(\xExt\), to the function space \(\Fint\), of functions defined on the interior grid \(\xInt\). Of course, \(\Fext = \RR^\Next\) and \(\Fint=\RR^\Nint\), but thinking of these spaces merely as high-dimensional Euclidean spaces neglects the meaning to which we have assigned the elements of each space. Further, we denote by \(\Fext_\vphi\) the set of functions in \(\Fext\) which, in addition, satisfy the discretized boundary conditions represented by the final \(2\abs{\EE}\) rows of system~\eqref{discretizedLaplacian}. Note that if \(\vphi\neq\Zero\), i.e., for nonhomogeneous vertex conditions, then \(\Fext_\vphi\) is an affine space of dimension \(\Nint\), while for \(\vphi=\Zero\), \(\Fext_\Zero\) is a linear vector space of dimension \(\Nint\). Thus the first \(\Nint\) rows of Eq.~\eqref{discretizedLaplacian} use the points from \(\xExt\) to approximately evaluate the underlying Laplace equations at the points in \(\xInt\), while the remaining \(2\abs{\EE}\) rows ensure that the solution lies on \(\Fext_\vphi\). 

In what follows, we will apply similar reasoning to discretize other problems on the quantum graph. As above, we apply the differential equations at the interior points and supplement these equations with \(2\abs{\EE}\) additional equations representing the vertex conditions, which suffice to specify a unique solution. In addition to the matrices \(\LVC\) and \(\PZero\), we will use the matrices
\begin{equation} \label{L0_PVC}
\LZero = \stackmat{\Lint} {\Zero_{2\abs*{\EE} \times \Next}}
\qand
\PVC = \stackmat{\Pint}{\MVC}.
\end{equation}
to develop time-stepping schemes for evolution equations in Sec.~\ref{sec:timestepping}.

Our choice of this ghost-point discretization stems from an effort to preserve the Laplacian's self-adjointness and, consequently, the realness of its eigenvalues post-discretization. Ghost points have an important advantage over the standard on-point discretization for enforcing the boundary conditions in the BVP~\eqref{robinSegment}. The simplest second-order discretization uses a non-staggered grid. It approximates the boundary condition at $x=0$ to second-order accuracy at the boundary by applying a one-sided difference to $\left.\dv{u}{x}\right\lvert_{x=0}$ using the values \(u_0\), \(u_1\), and \(u_2\), and similarly at the $x=1$. Solving the two discretized boundary conditions eliminates the values of \(u_0\) and \(u_{N+1}\) from the system, leaving a system of \(N\) unknowns, but with an asymmetric finite-difference matrix, even though the differential operator which it approximates is, self-adjoint. The ghost-point discretization, by contrast, preserves self-adjointness (after the two ghost points are first eliminated from the system). 

This changes slightly when we extend this construction to the quantum Graph Poisson problem~\eqref{Poisson}.
If all edges are discretized with the same step size \(h\), the discretization matrix constructed above is symmetric. However, choosing all discretization lengths to be equal may be impossible or inconvenient. In that case, we may measure the magnitude of the asymmetry by considering the largest element, in absolute value, of the asymmetric part \(\frac{1}{2} \left( \A - \transpose{\A} \right) \). Assume that edge \(\edge_m\) is discretized with a step size \(h_m = h + \delta_m\) with \(\delta_m = O(\delta) \ll h\). Then using one-sided centered differences introduces terms of \(O\left(\frac{1}{h^2}\right)\) into the asymmetric part of \(\A\). By contrast, using ghost points introduces terms of \(O\left(\frac{\delta}{h^2}\right)\). Therefore, if all the discretization lengths \(h_m\) are roughly equal, the non-symmetric part of \(\A\) will be significantly smaller so that the matrix \(\A\) is "more symmetric." Since \(\A\) is not symmetric, we cannot guarantee that all of our eigenvalues are purely real, as they are for the underlying differential operator. However, by an informed choice of the discretization sizes, we may minimize the effects of the asymmetry

\paragraph{\textbf{An example discretization}}

\begin{figure}[t]
\begin{center}
\begin{subfigure}[b]{0.14\textwidth}
    \includegraphics[width=\textwidth]{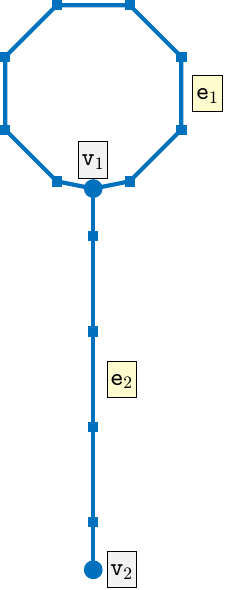}
    \caption{}
\end{subfigure}  \quad  \quad  \quad
\begin{subfigure}[b]{0.42\textwidth}
    \includegraphics[width=\textwidth]{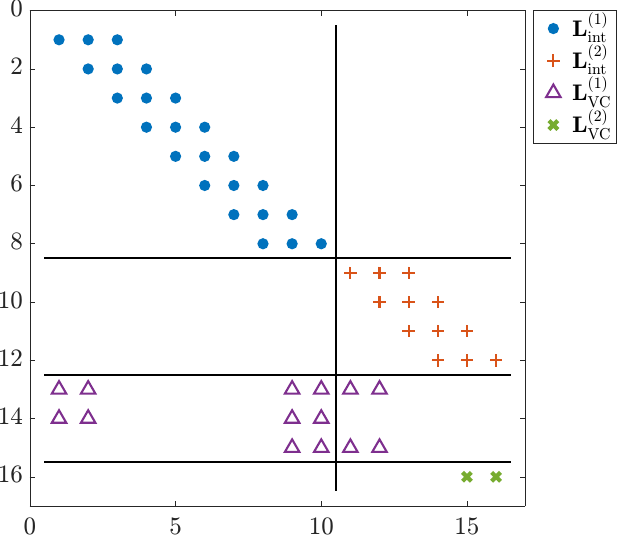}
    \caption{}
\end{subfigure}
\end{center}
   \caption{(a) The lasso graph, showing interior discretization points. (b) The structure of the nonzero entries in the matrix \(\LVC\), the Laplacian matrix extended with vertex conditions.}
   \label{fig:lasso}
\end{figure}

We display the structure of the discretized Laplacian matrix of a lasso graph with two vertices and two edges in Fig.~\ref{fig:lasso}. The edge \(\edge_1\) points from \(\vertex_1\) to \(\vertex_2\) and the edge \(\edge_2\) points from \(\vertex_2\) to itself. The discretization has \(N_1=4\) and \(N_2=8\)  points. The figure shows the interior points of the discretization and the vertices but not the ghost points. The figure also shows the nonzero structure of the matrix \(\LVC\). The \(16\times2\) matrix \(\MNH\)  is nonzero at \((13,1)\) and \((14,2)\).

\subsubsection{Chebyshev discretization}

To achieve spectral accuracy, QGLAB  allows discretization using \emph{rectangular collocation}, a method based on Chebyshev polynomials due to Driscoll and Hale~\cite{driscoll2016rectangular}, with further implementation details described in Ref.~\cite{Xu:2016dz}. Enforcing non-trivial boundary conditions based on a non-square differentiation matrix with this method is especially simple. To introduce the idea, we again consider the Robin problem on a line segment defined in Eq.~\eqref{robinSegment}. The exterior grid \(\xExt\) is given by the \(N+2\) (increasing) Chebyshev points of the second kind,
\begin{equation}
    \xkExt{k} = \frac{\ell}{2}\left(1-\cos{\left(\frac{k\pi}{N+1}\right)}\right), \qquad k=0,1,\ldots,N, N+1.
    \label{chebpts1}
\end{equation}

We adapt the notation of Eq.~\eqref{u_f_phi} to define the vectors \(\u\), \(\f\), and \(\vphi\) on the discretization points defined in Eq.~\eqref{chebpts1}. The main observation motivating rectangular collocation is that applying a second derivative matrix defined over a finite space of polynomials should reduce the order of that space by two, naturally leading to matrices of size \(N\times(N+2)\). This is realized by first operating on the vector \(\u\) with the standard \((N+2)\times(N+2)\) Chebyshev derivative matrix \(\DD\) and then resampling these polynomials onto the interior grid \(\xInt\)   of first-kind Chebyshev points
\begin{equation}
    \xkInt{k}=\frac{\ell}{2}\left(1- \cos{\left( \frac{(2k-1)\pi}{2N} \right)} \right) \qquad k=1,2,\ldots,N.
\end{equation}

As in the finite-difference discretization above, we  define exterior and interior grids \(\xExt\) and \(\xInt\). The approximate solutions are defined on \(\xExt\), but derivatives, and thus approximations to the differential equations, are evaluated at the points of \(\xInt\). By contrast, the two grids are disjoint sets and the first and last points of the exterior grid \(\xExt\) are the endpoints of the interval, not the ghost points.

\begin{figure}[htbp]
    \centering
    \includegraphics[width=0.9\textwidth]{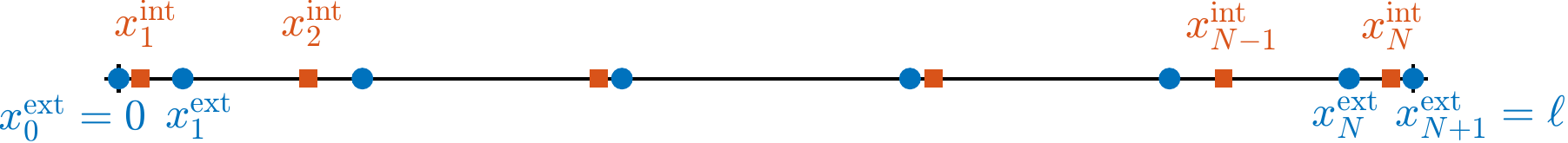}
    \caption{Discretization of the interval \([0,\ell]\) using a grid \(\xExt\) of Chebyshev points of the second kind (in blue) and a grid \(\xInt\) first kind Chebyshev points (in red).}
    \label{fig:chebpts}
\end{figure}

Resampling is a linear operation and is represented by an \(N\times(N+2)\)-dimensional \emph{barycentric resampling matrix} \(\Pint\) whose construction uses the barycentric interpolation formula proposed in~\cite{berrut2004barycentric}. Given the set of points \(\xExt = \{\xkExt{k}\}_{k=0}^{N+1}\), the barycentric weights are
\begin{equation}
    {\tilde w}_k = \prod_{\substack{l=0\\l\neq k}}^{N+1} (\xkExt{k} + \xkExt{l})^{-1}, \qquad k=0,\dots,N+1.
\end{equation}
These are used to construct a unique interpolating polynomial \(p_{N+1}(x)\) which interpolates the set of data points \(\{(\xkExt{k},f_k)\}_{k=0}^{N+1}\). The polynomial is evaluated at both \(\{\xkExt{k}\}_{k=0}^{N+1}\) and \(\{\xkInt{k}\}_{k=1}^{N}\) so that the barycentric resampling matrix is given by
\begin{equation}
    (\Pint)_{j,k} = 
    \begin{cases}
    \frac{{\tilde w}_k}{\xkInt{j} - \xkExt{k}}\left( \sum_{l=0}^{N+1}\frac{{\tilde w}_l}{\xkInt{j}-\xkExt{l}} \right)^{-1} &  \xkInt{j} \neq \xkExt{k},\\
     1 & \xkInt{j} = \xkExt{k},
     \end{cases}
\end{equation}
and satisfies 
\begin{equation}
    p_{N+1}(\xInt) = \Pint~p_{N+1}(\xExt),
\end{equation}
i.e., the barycentric resampling matrix maps a polynomial's values at the gridpoints $\xExt$ to its values at the gridpoints $\xInt$.  Driscoll and Hale generalize this in~\cite{driscoll2016rectangular}.

Putting this all together, the product  
\begin{equation}\label{LintCheb}
    \Lint = \Pint \DD
\end{equation}
defines a \(N\times (N+2)\) differentiation matrix. The right-hand side of the differential equation~\eqref{robinSegment} must be resampled to the same grid, so the differential equation is discretized by the \(N\) equations \(\Lint \u = \Pint \vphi\), leaving two equations to define the boundary conditions~\eqref{robinSegment}. These may be compactly rewritten as
\[
\MBC \u = \begin{pmatrix} \phi_0 \\ \phi_L \end{pmatrix},
\]
where \(\MBC\) is a matrix of size \(2\times(N+2)\) conveniently expressed using unit vectors
\begin{equation}
\MBC = \def\arraystretch{1.25}
\begin{pmatrix}
\e_1^\Tee \D + \alpha_0 \e_1^\Tee \\
-\e_{N+2}^\Tee \D + \alpha_L \e_{N+2}^\Tee
\end{pmatrix}.
\label{MBC_cheb}
\end{equation}
The matrices \(\Lint\) and \(\Pint\) are blockwise dense, in contrast to the banded matrices arising in the uniform discretization, defined in Eq.~\eqref{MBC}. We have constructed all the necessary elements to reinterpret Eq.~\eqref{discSystemEqs} as a spectral collocation of Eq.~\eqref{robinSegment}.

\begin{figure}[h] 
   \centering
\begin{subfigure}[b]{0.13 \textwidth}
    \includegraphics[width=\textwidth]{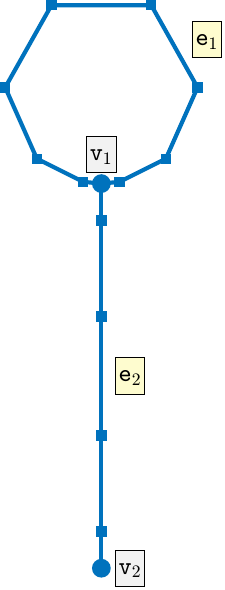}
    \caption{}
\end{subfigure} \quad \quad
\begin{subfigure}[b]{0.4\textwidth}
    \includegraphics[width=\textwidth]{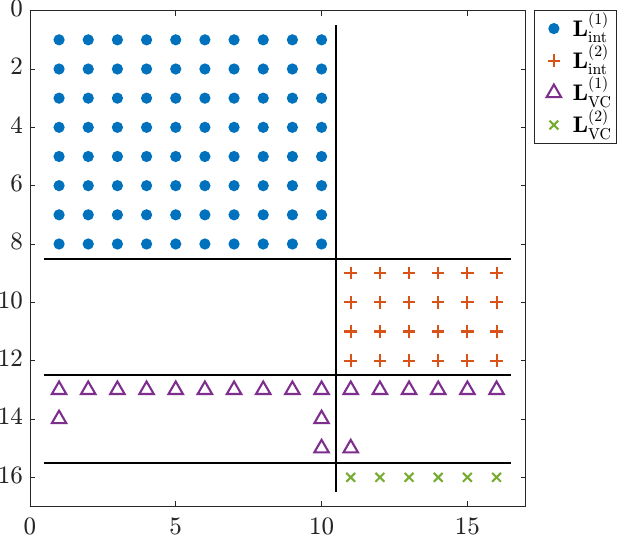}
    \caption{}
\end{subfigure}
\begin{subfigure}[b]{0.4\textwidth}
    \includegraphics[width=\textwidth]{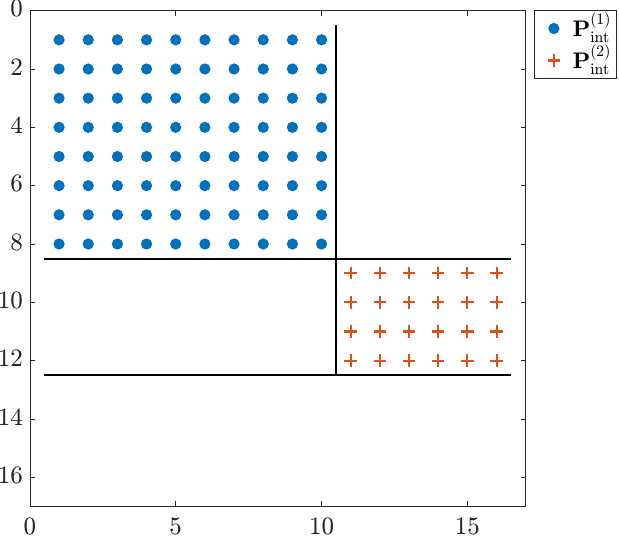}
    \caption{}
\end{subfigure}
   \caption{(a) The lasso graph, shown with interior discretization points in the Chebyshev discretization. (b) The structure of the nonzero entries in \(\LVC\), the Laplacian matrix extended with vertex conditions. (c) The structure of \(\P_0\), the barycentric resampling matrix extended with zeros.}
   \label{fig:lassocheb}
\end{figure}

Extending this argument from the ODE boundary problem proceeds as for the centered-difference approximation. 
Defining the matrices \(\LVC\), \(\PZero\), and \(\MNH\) in Eq.~\eqref{discretizedLaplacian} is straightforward once we construct the submatrix \(\MVC^{(n)}\) defining the discretized vertex condition~\eqref{Luniform} and \eqref{Puniform}, extending the construction in~\eqref{MBC_cheb}.  Fig.~\ref{fig:lassocheb} shows a coarse discretization of the example shown in Fig.~\ref{fig:lasso}, in which the blocks defining the second derivative are dense, as are the rows defining the Robin-Kirchhoff vertex condition, whereas the rows enforcing continuity contain only two nonzero entries.

Two of the authors used this method to solve the Laplacian eigenproblem on an interval perturbed by several delta function potentials in~\cite{beck2020limiting}.

\subsection{Numerical and Symbolic Eigenproblems}
\label{sec:numerical_eigenproblems}
\begin{figure}[htbp] 
   \centering
   \includegraphics[width=\textwidth]{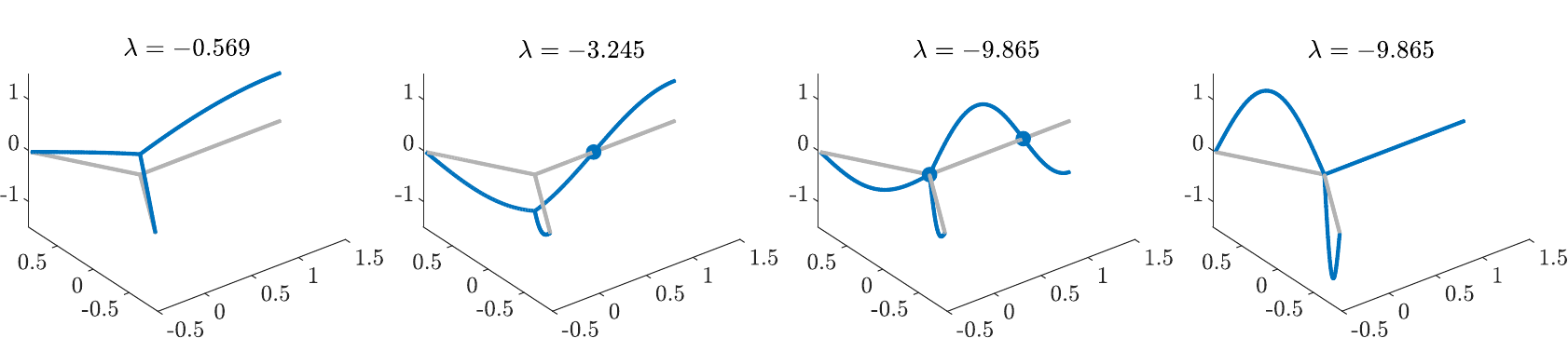} 
   \caption{Four eigenfunctions of a Y-shaped quantum graph.}
   \label{fig:star_eigenfunctions}
\end{figure}

Following the steps used above to discretize the Poisson problem,  the  eigenvalue problem~\eqref{analyticEvalProb} becomes
\begin{equation}\label{numericEvalProb}
    \LVC \u = \lambda \PZero \u
\end{equation}
in both the finite-difference and Chebyshev discretizations. Since \(\PZero\) is not the identity matrix and is singular, this is a \emph{generalized eigenvalue problem} which can be solved using MATLAB's \eigs.
QGLAB has overloaded the \eigs\ command so that \listtt{[d,v]=G.eigs(m)} returns the \(m\) eigenvalues of the smallest absolute value. Fig.~\ref{fig:star_eigenfunctions} shows the eigenfunctions of the four smallest eigenvalues of the Laplacian on a Y-shaped graph with Dirichlet conditions at the ends of the two shorter edges. The example is described further in Sec.~\ref{sec:eigen_example}. QGLAB's plotting features are described in Sec.~\ref{sec:functions-and-plotting}.

The zeros \(k_n\) of the secular determinant function \(\Sigma(k)\), described in Sec.~\ref{sec:eigenvalueProblem},  correspond to eigenvalues \(\lambda_n=-k_n^2\) of the Laplacian operator.  The QGLAB function \listtt{secularDet}  uses MATLAB'S Symbolic Math Toolbox to construct \(\Sigma(k)\) for the boundary conditions~\eqref{KirchhoffRobin} with \(w_j \equiv 1\) and for Dirichlet boundary conditions and simplify it for typesetting and plotting. This is described in the dissertation~\cite{Conte:2022}. It provides a check on the numerical calculation of eigenvalues for the discretized problem.

\subsection{Nonlinear Solvers, Continuation, and Bifurcation Algorithms}
\label{sec:solvers}

\begin{figure}[tbp] 
   \centering
      \begin{subfigure}[b]{0.3\textwidth}
          \includegraphics[width=\textwidth]{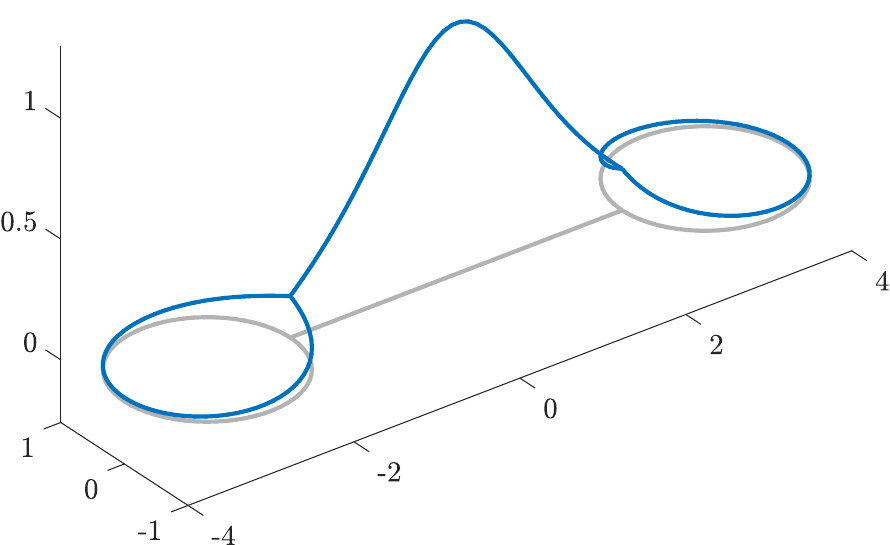} 
          \caption{}
      \end{subfigure}
\begin{subfigure}[b]{0.3\textwidth}
     \includegraphics[width=\textwidth]{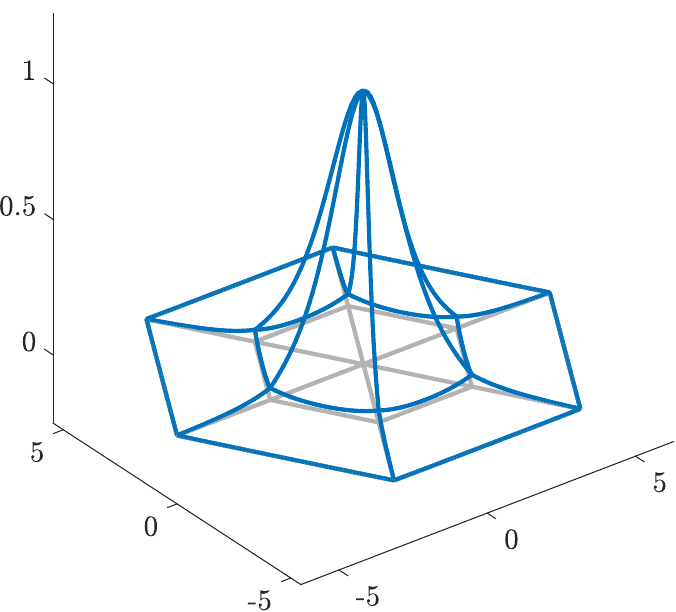}
     \caption{}
\end{subfigure}
   \caption{(a) A standing wave of cubic NLS on a dumbbell graph with \(\Lambda=-1\). (b) A standing wave on a spiderweb graph with \(\Lambda=-1\).}
   \label{fig:stationaryNLS}
\end{figure}

After discretizing the spatial derivatives, we use QGLAB to construct the Newton-Raphson method to find standing wave solutions of the stationary NLS~\eqref{stationary_NLS} with $\sigma=1$. Fixing $\Lambda$, it proceeds by the iteration $\Psi_{n+1} = \Psi_{n} + \delta$ where $\delta \in \Fext_0$ solves $D\FF(\Psi_n,\Lambda)\cdot \delta = -\FF(\Psi_n,\Lambda)$. Enforcing the vertex conditions gives an iteration of the form
\begin{equation}\label{Newton}
\left( \LVC + \PZero \left(6 \diag(\Psi_n^2) + \Lambda \I \right) \right) \delta = 
- \LZero \Psi_n - \PZero \left(\Lambda \Psi_n + 2 \diag(\Psi_n^3) \right).
\end{equation}
Two examples of such solutions are shown in Fig.~\ref{fig:stationaryNLS} and Sec.~\ref{sec:NLS_example}

Solutions to Eq.~\eqref{stationary_NLS} do not occur at isolated points but along one parameter families that, away from singularities, can be parameterized by the frequency \(\Lambda\). Pseudo-arclength continuation provides a way to follow this family as it traces a smooth path. It is due originally to Keller~\cite{Keller:1977} and is well summarized in the textbook of Nayfeh and Balachandran~\cite{nayfeh:1995}, who cite many additional contributors to the method. The method allows curves to be continued around fold singularities, and other techniques can detect branch points, which encompass both pitchfork and transcritical bifurcations, and initiate new families of solutions branching off from the computed branch; see also Govaerts~\cite{Govaerts:2000vv}.  QGLAB's pseudo-arclength continuation was first used to calculate bifurcation diagrams in Ref.~\cite{Berkolaiko:2020}.

QGLAB's contribution here is not a novel numerical method but the integration of pseudo-arclength continuations into quantum graph software. Failure to use these methods has led to some seemingly impossible phenomena in the quantum graph literature: a solution branch that appears to end abruptly in Ref.~\cite{Marzuola:2016bl}, subsequently fixed in Ref.~\cite{Goodman:2019go}, and solutions with different symmetries seeming to lie on a single branch in Ref.~\cite{Besse:2021}, fixed in the example below.

The main step of the continuation algorithm is a Newton step like~\eqref{Newton}, in which the parameter \(\Lambda\) is left unknown. A predictor-corrector algorithm introduces one more equation to close the system. QGLAB is novel among quantum graph software because it implements continuation algorithms and because the operator language developed in Sec.~\ref{sec:discretization} allows a compact representation of the equations that arise.  

\begin{figure}[htbp] 
   \centering
         \begin{subfigure}[b]{0.3\textwidth}
          \includegraphics[width=\textwidth]{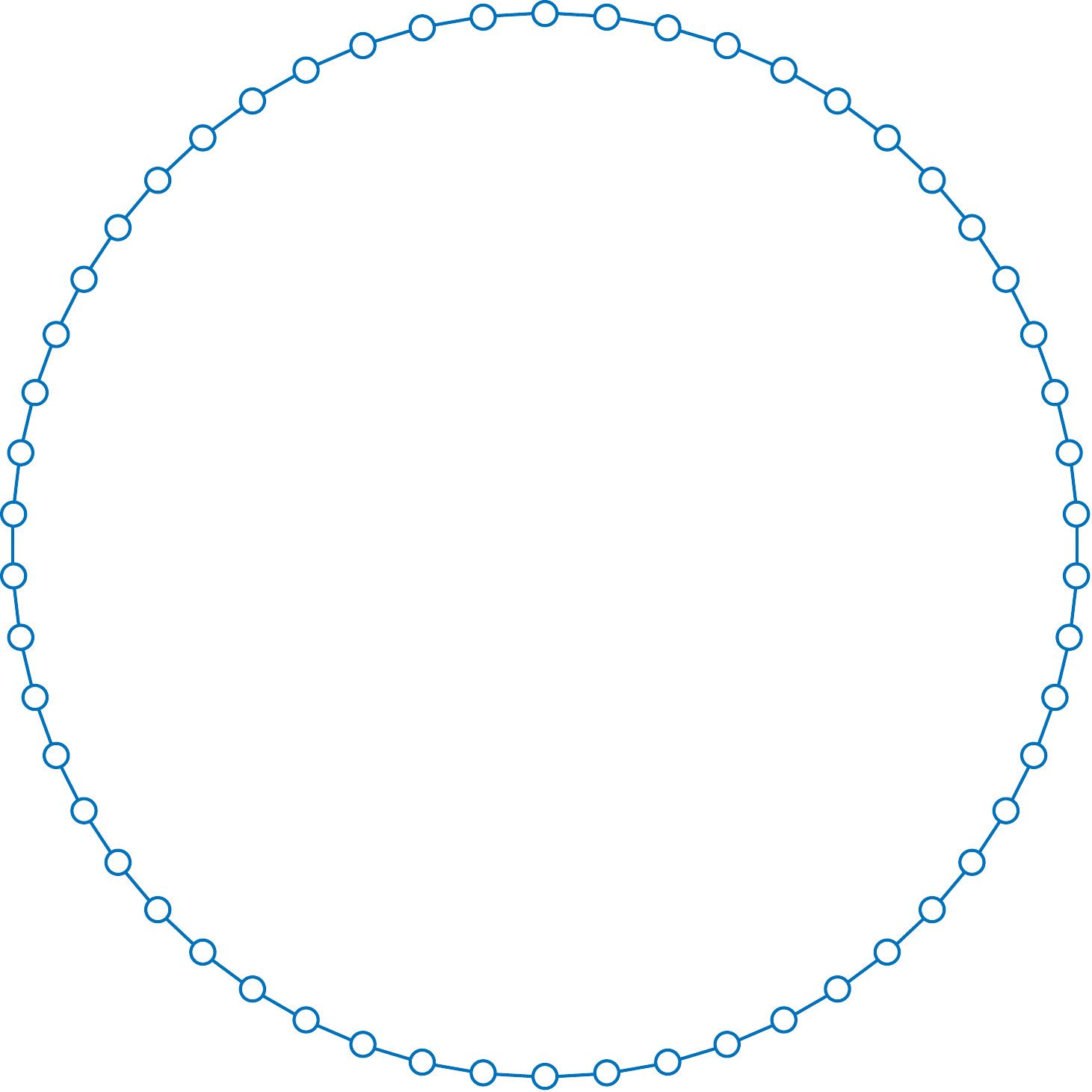} 
          \caption{}
      \end{subfigure}
\begin{subfigure}[b]{0.3\textwidth}
     \includegraphics[width=\textwidth]{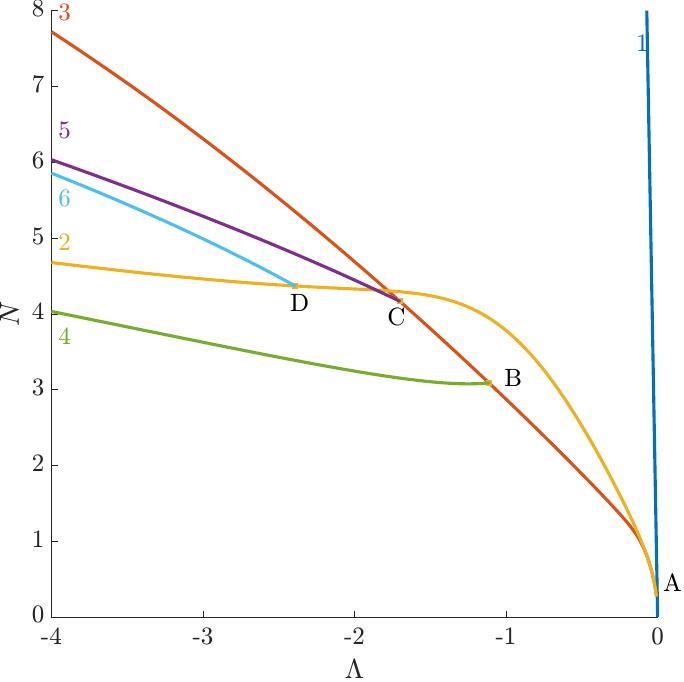} 
     \caption{} 
\end{subfigure} 
\begin{subfigure}[b]{0.9\textwidth}
      \includegraphics[width=\textwidth]{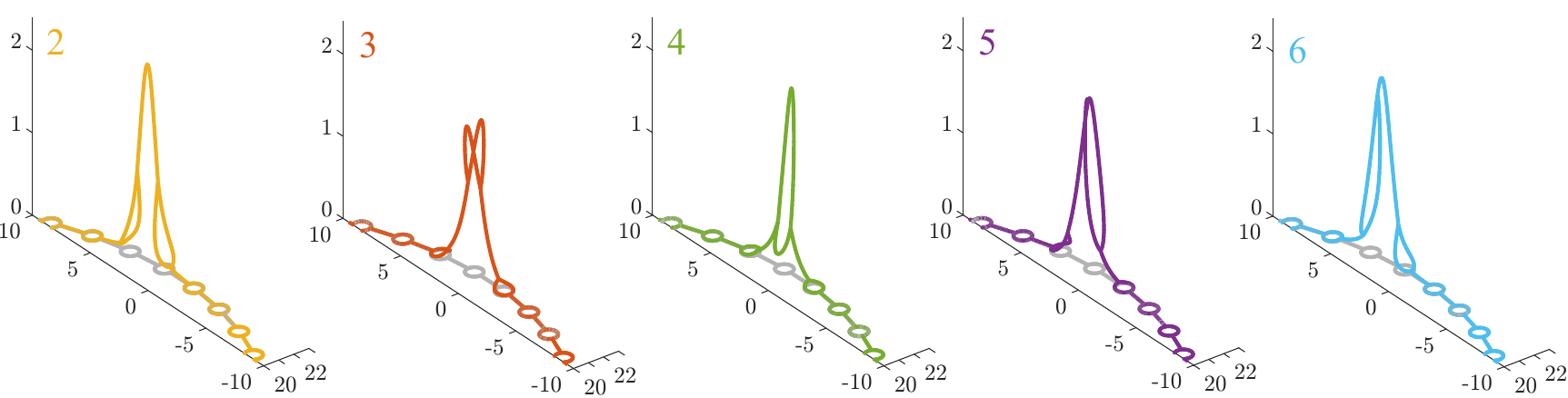} 
      \caption{}
\end{subfigure}
   \caption{(a) Layout of the necklace graph, with 54 "strings" and 54 "pearls." (b) Partial bifurcation diagram plotting the solutoin's mass as a function of its frequency. (c) Solutions along the color-coded branches with frequency $\Lambda \approx -4$.}
   \label{fig:necklace}
\end{figure}

An example computation demonstrates QGLAB's capabilities. Sec.~\ref{sec:nonlinear_standing_examples} of the \suppendix\ contains further examples, but these computations involve too many lines of code to include in the published article, so these programs are understood most easily from the live script \listtt{continuationInstructions.mlx}.
Fig~\ref{fig:necklace}(a) shows a so-called necklace graph, similar to an example in Ref.~\cite{Besse:2021}, consisting of 54 segments alternating between "strings" and "pearls." The cited paper allows the numerical domain to widen as the amplitude decreases, demonstrating the ground state scaling at small amplitude. However, because it does not employ continuation and branch-switching, it does not demonstrate the relationships between the branches.

Panel (b) of the figure shows a partial bifurcation diagram of standing wave solutions to cubic NLS, plotting $N = \norm{\Psi}_2^2$ as a function of $\Lambda$. Branch \textbf{1} contains solutions with $\Psi(x,\Lambda)$ constant on $\Gamma$. At the point marked \textbf{A}, branches \textbf{2} and \textbf{3} bifurcate from the first branch. Solutions on branch \textbf{2} are symmetric around their maxima on the center of a string and those on branch \textbf{3} are symmetric around the centers of the pearls. At points \textbf{B}, \textbf{C}, and \textbf{D}, new branches arise due to symmetry-breaking (pitchfork) bifurcations. The figure shows five solutions with frequency $\Lambda \approx -4$. Solution \textbf{4} is the minimal mass solution at that frequency.

\subsection{Time-stepping for evolution problems posed on a quantum graph}
\label{sec:timestepping}

Our goal in this section is to adapt well-known time-dependent PDE solvers to work on quantum graphs using the framework described above for stationary problems. The first subsection contains two examples, while the second discusses the construction of a general solver for a certain class of such problems.

\subsubsection{Elementary methods}
\label{sec:elementary_time}

We first adapt two standard methods to quantum graphs: the Crank-Nicholson method for the heat equation and the leapfrog method for the sine-Gordon equation. In the following, $\tau$ is the time step, $t_n = n\tau$ is the discretized times, $\psi_n$ represents the solution at time $n$ discretized in time only and $\vpsi_n$ the spatially discretized solution at time $n$.

\paragraph{\textbf{Crank-Nicholson for the heat equation}}
The Crank-Nicholson method is a common second-order method for the heat equation
\(
\pdv{\psi}{t} =  \triangle \psi. 
\)
Discretizing only in time, the update at time $t_{n+1}$ is found by solving the equation
\[
\frac{\psi_{n+1} - \psi_{n}}{\tau} =  
  \frac{ \triangle  \psi_{n} + \triangle  \psi_{n+1}}{2}.
\]
 In QGLAB, this is discretized and evaluated on the interior grid to give
\[
\left(\Pint - \frac{\tau}{2} \Lint \right) \vpsi_{n+1} = 
\left(\Pint + \frac{\tau}{2} \Lint \right) \vpsi_{n}.
\]
Combining this with homogeneous vertex conditions yields
\begin{equation} \label{CNstepper}
\Lminus \vpsi_{n+1} \equiv
\left(\PVC - \frac{\tau}{2} \LZero \right)
 \vpsi_{n+1} = 
\left(\PZero + \frac{\tau}{2} \LZero \right) \vpsi_{n}
= \Lplus \vpsi_n.
\end{equation}

The method iterates the MATLAB code \listtt{y = Lminus \\ (Lplus*y)}. An example in Sec.~\ref{sec:heat_example} in the \suppendix\ computes a solution to the heat equation on the dumbbell graph and demonstrates convergence.

\paragraph{\textbf{Leapfrog for nonlinear Klein-Gordon equations}}

The following example, contained in the live script \listtt{sineGordonOnTetra.mlx}, solves the sine-Gordon equation on the tetrahedron quantum graph, considered previously in~\cite{Dutykh:2018}, which consists of a wave equation with a sinusoidal nonlinearity,
\begin{equation}
\psi_{tt} -\triangle \psi + \sin{\psi} = 0.
\label{sineGordon}
\end{equation}
Discretizing in time only and applying second-order centered differences in time gives 
\begin{equation}
\frac{\psi_{n+1} - 2 \psi_n +  \psi_{n-1} }{\tau^2}= \triangle \psi_n - \sin{\psi_n}.
\label{leapfrog}
\end{equation}
Applying the discretization in space and solving gives
\begin{equation*} 
\Pint \vpsi_{n+1} = \Pint \left( 2\vpsi_n - \vpsi_{n-1} - \tau^2 \sin{\vpsi_n} \right) + \tau^2 \Lint \vpsi_n 
\end{equation*}
and enforcing the vertex conditions gives
\begin{equation}\label{leapfrog_step}
\PVC \vpsi_{n+1} = \PZero \left( 2\vpsi_n - \vpsi_{n-1} - \tau^2 \sin{\vpsi_n} \right) + \tau^2 \LZero \vpsi_n.
\end{equation}
The iteration requires an approximation $\vpsi_1 \in \Fext_0$ at time $t_1$, which may be found from the initial conditions
 \(\left.\vpsi\right\rvert_{t=0} = \vpsi_0\) and \(\left.\pdv{t} \vpsi\right\rvert_{t=0} = \vphi_0\)  using \(O(\tau^2)\) approximation
\[
\PVC \vpsi_1  = \PZero \left( \vpsi_0 + \tau \vphi_0 - \frac{\tau^2}{2} \sin{\vpsi_0} \right)
  + \frac{\tau^2}{2} \LZero \vpsi_0.
\]
The time-stepper in Eq.~\eqref{leapfrog_step} is implemented by the line
\begin{lstlisting}[numbers=none]{matlab}
u2 = PVC \ (P0*(2* u1 - u0 - tau^2*sin(u1)) + tau^2*L0*u1).
\end{lstlisting}

In Sec~\ref{sec:sG_example} in the \suppendix, we compute an example from Section $4.1$ of Ref.~\cite{Dutykh:2018} in which solitary waves collide with vertices on a tetrahedron metric graph and demonstrate convergence.

\subsubsection{A general purpose higher-order time stepper}
\label{sec:timesteppers}

We now construct a general-purpose solver for differential equations of the form
\begin{equation} 
\pdv{\psi}{t} =  \mu \Laplacian \psi + f(\psi),
\label{evolution}
\end{equation}
posed on the quantum graph subject to any homogeneous vertex conditions implemented in QGLAB and such that \(f(\psi)\) contains any terms involving a potential or nonlinearity. Depending on the constant \(\mu\), which could be real, imaginary, or complex, this formulation includes heat, Schrödinger, Ginzburg-Landau, and scalar reaction-diffusion equations. The vertex condition constraints make the straightforward application of standard methods somewhat difficult. Here, we construct a method that overcomes these problems.

The main issues in constructing a time-stepping algorithm for an evolutionary PDE defined on a quantum graph using the spatial discretization described in Sec.~\ref{sec:discretization} can be illustrated using Euler methods.  These ideas then extend straightforwardly to Runge-Kutta algorithms. 
The forward Euler method is 
\[
\psi_{n+1}=\psi_{n}+\tau  \left(  \mu \Laplacian \psi_{n}+f( \psi_{n})\right),
\]
subject to vertex conditions applied to \(\psi_{n+1}\).
After we discretize in space and enforce the vertex conditions, this yields a time-stepper
\begin{equation} \label{forward_euler}
\PVC \vpsi_{n+1}=\PZero \cdot \left( \vpsi_{n}+\tau  \f\left(\vpsi_{n}\right) \right) +\tau  \mu \LZero \vpsi_{n}.
\end{equation}
The matrix \(\PVC\) on the left makes the method implicit, but the implicitness is linear, requiring no Newton iterations. However, the stiff Laplacian term on the right is evaluated explicitly, imposing a step-size restriction, so the method is impractical.

Similarly, the backward Euler method is
\[
\psi_{n+1}=\psi_{n}+\tau  \left(  \mu \Laplacian \psi_{n+1}+f(\psi_{n+1})\right),
\]
subject to vertex conditions on \(\psi_{n+1}\).
After we discretize in space and enforce the vertex conditions, this yields a time-stepper
\begin{equation} \label{backward_euler}
\left(\PVC-\tau  \mu \LVC\right) \vpsi_{n+1}-\tau  \PZero \f(\vpsi_{n+1}) =\PZero \vpsi_{n}.
\end{equation}
The stiff laplacian term is handled implicitly and imposes no time-step restrictions, but the implicit nonlinear term requires Newton iterations at each step and slows down the method.

To resolve this difficulty, we may treat the stiff term involving the Laplacian implicitly and the nonstiff term involving the nonlinearity explicitly. This idea was introduced for Runge-Kutta methods, which include the Euler method, by Ascher et al.~\cite{Ascher:1997}. There exist several such implicit-explicit (IMEX) Euler methods, including one they call forward-backward Euler \((1,1,1)\):
\[
\psi_{n+1}=\psi_{n}+\tau  \left(  \mu \Laplacian \psi_{n+1}+f( \psi_{n})\right)
\]
subject to vertex conditions on \(\psi_{n+1}\).
After we discretize in space and enforce the vertex conditions, this yields a time-stepper
\begin{equation} \label{imex_euler}
\left(\PVC  -\tau  \mu \LVC \right) \vpsi_{n+1}=\PZero \cdot \left(  \vpsi_{n}+\tau  \f(\vpsi_{n})\right).
\end{equation}
This method combines the best aspects of the forward and backward Euler methods. Moving the operator \(\tau \mu \LVC\) to the left-hand side resolves the stiffness issue without requiring \(\tau\) to be small. Keeping the nonlinear term on the right-hand side eliminates the need to solve a nonlinear equation on each step. However, the method is only first order in time, requiring small time steps for accuracy. 

Ref.~\cite{Ascher:1997} applies similar ideas to derive implicit-explicit (IMEX) Runge-Kutta methods, which at each stage handle the stiff part of the evolution equation implicitly and the nonstiff part explicitly. The implicit terms are \emph{strongly diagonal}, so each substage of a time step can be found in terms of the previous substages. QGLAB comes with a four-stage third-order Runge Kutta method \listtt{qgdeSDIRK443}, based on the method denoted \((4,4,3)\) in~\cite{Ascher:1997}.

\begin{figure}[htbp] 
   \centering
   \includegraphics[width=0.8\textwidth]{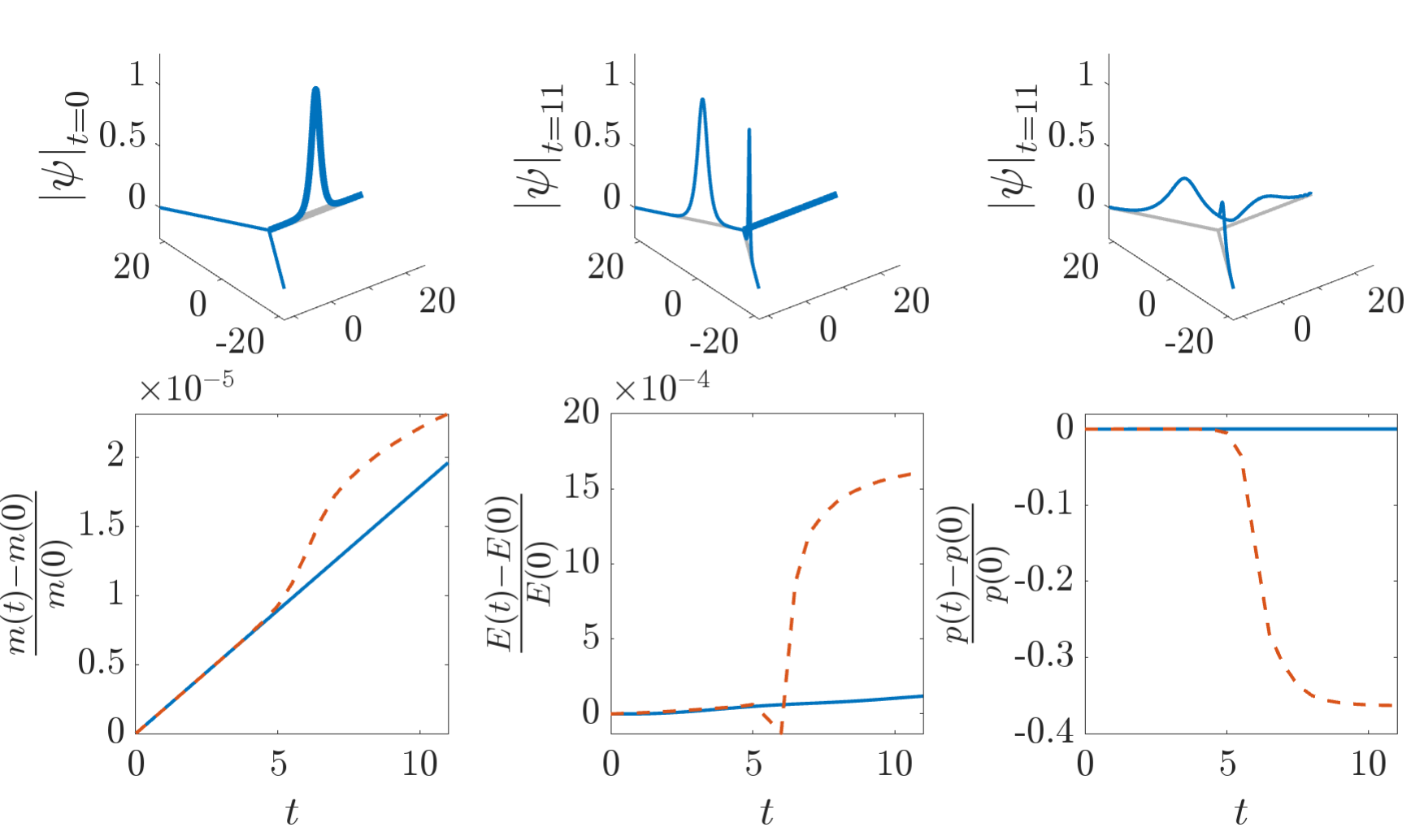} 
   \caption{Collision of an NLS soliton with the vertex of a star graph.  Top: (Left) Initial time. (Center) Final time, balanced graph. (Right) Final time, unbalanced graph. Bottom: Relative change in mass (Left), energy (Center), and momentum (Right) with balanced graph in solid blue and the unbalanced graph in dashed red.}
   \label{fig:NLSonStar}
\end{figure}

Figure~\ref{fig:NLSonStar} shows the collision of a soliton with a vertex on a star graph with both so-called balanced and unbalanced Kirchhoff conditions, as considered in Ref.~\cite{Kairzhan:2019fe} and explained further in Sec.~\ref{sec:NLS_example}. The soliton passes through the balanced vertex but is largely reflected by the unbalanced vertex. NLS dynamics conserves the mass (squared $L^2$ norm~\eqref{Lpnorm}) and energy~\eqref{energyNLS}; for certain initial conditions, the dynamics on the balanced graph also conserve momentum.

\section{Understanding the MATLAB implementation}
\label{sec:MATLAB}

\subsection{MATLAB's digraph class}
\label{sec:digraph}
MATLAB provides tools for computations on undirected and directed graphs. The main component of QGLAB is a class \qGraph, which builds upon MATLAB's \digraph\ class used for defining directed graph objects. The graph in Fig.~\ref{fig:sillyGraph} is constructed using:
\begin{lstlisting}[numbers=none]{matlab}
source=[1 1 1 2 2]; target=[1 1 2 2 3];
G = digraph(source, target);
plot(G)  % Actual figure created using quantumGraph plot function
\end{lstlisting}
The vectors \source\ and \target\ specify, respectively, the initial and final vertices of the four edges. The digraph object \ttG\ contains two fields: \listtt{G.Edges} and \listtt{G.Nodes}, each in the form of a table---a MATLAB array type that holds column-oriented data, each column stored as a variable. The methods \listtt{G.numnodes} and \listtt{G.numedges} return the number of vertices (nodes) and edges, respectively. The array of edges contains one variable \listtt{G.Edges.EndNodes}, an \( \abs*{\EE} \times 2 \) array whose two columns contain, respectively, the indices of the source vertices and the target vertices. The nodes table is initially empty. We add fields to the two tables to create the quantum graph class.

\subsection{Understanding the quantum graph class and initializing a quantum graph object}
\label{sec:initializing}
The graph and matrix shown in Fig.~\ref{fig:lasso} were generated with the code:
\begin{lstlisting}[numbers=none]{matlab}
source=[1 2]; target=[2 2]; L=[4 2*pi]; nx=[4 8];
G=quantumGraph(source,target,L,'nxVec',nx)
\end{lstlisting}
The last line initializes a \qGraph\ object. The three required arguments \source, \target, and \ttL\ must be entered in that order. The last, \ttnx, is an optional argument. If \ttnx\ is a vector of length \listtt{G.numedges}, it defines the number of interior points on each edge. If it is a scalar, the constructor will assign \ttnx\ points per unit length to each edge, rounding if necessary. 

The constructor takes several optional arguments, which will be discussed below. Some have default values if not specified in the function call. Optional arguments are listed in the function call using a key/value syntax after required arguments. In older releases of MATLAB, this is entered as 
\lstinline|G=quantumGraph(source,target,L,'key1',value1,'key2',value2)|,
while more recent releases allow the more compact syntax
\lstinline|G=quantumGraph(source,target,L,key1=value1,key2=value2)|.
Complete instructions, including optional arguments, are presented in Sec.~\ref{sec:function_listing} of the \suppendix.

The above commands return the following in the MATLAB command window:
\begin{lstlisting}[numbers=none]{matlab}
G =  quantumGraph with properties:
  
            discretization: 'Uniform'
       wideLaplacianMatrix: [12x16 double]
       interpolationMatrix: [12x16 double]
          discreteVCMatrix: [4x16 double]
    nonhomogeneousVCMatrix: [16x2 double]
          derivativeMatrix: [16x16 double]
                 potential: []          
\end{lstlisting}
The most important property of this \qGraph\ object is \qg, which specifies the quantum graph itself and its discretization. It is not visible in the above listing because it is a \emph{private} property of the object. The user can not directly access it, but class methods may act on it. We will discuss it last. The remaining properties are publicly viewable but can only be set by class methods. They are:
\begin{description}
\itemtt{discretization} may take three values, \Uniform (default), \Chebyshev, or \none. With the \none\ property, only secular determinant computations function.
\itemtt{wideLaplacianMatrix} The rectangular Laplacian matrix \(\Lint\) defined in Eq.~\eqref{Luniform} for either the uniform or Chebyshev discretization.
\itemtt{interpolationMatrix} The interpolation or resampling matrix \(\Pint\) as defined in Eq.~\eqref{Puniform} for either the uniform or Chebyshev discretization.
\itemtt{discreteVCMatrix} The matrix \(\MVC\) defining the discretization of the vertex conditions defined in Eq.~\eqref{Luniform} and \eqref{Puniform} for either discretization.
\itemtt{nonhomogeneousVCMatrix} The matrix \(\MNH\) defined in Eq.~\eqref{MNH} which maps nonhomogeneous terms in the vertex condition to the appropriate row.
\itemtt{derivativeMatrix} The square matrix used to calculate the first derivative on each edge. It is used to compute time-dependent solutions' energy and momentum functionals and to plot solution branches in continuation problems.
\itemtt{potential} The optional potential $V(x)$.
\end{description}

The property \qg\ is a MATLAB directed-graph object consisting of a \listtt{Nodes}\ table and an \listtt{Edges}\ table with added fields needed to define a quantum graph. 
Because \qg\ is a private property, viewing these tables using the syntax \listtt{G.qg.Nodes} and \listtt{G.qg.Edges} is disabled. Instead, \listtt{Nodes}\ and \listtt{Edges}\ \listtt{quantumGraph} methods have been written that return each of these tables, so we may view the tables using the syntax \listtt{G.Nodes} and \listtt{G.Edges}, as in this code listing. In addition, a method exists that returns each default table column; for example, the Robin coefficients can be returned by \listtt{G.robinCoeff}. We examine the node data, which has two fields
\begin{lstlisting}[numbers=none]{matlab}
>> disp(G.Nodes)
    robinCoeff     y 
    __________    ___
        0         NaN
        0         NaN
\end{lstlisting}
The fields are: 
\begin{description}
  \itemtt{robinCoeff} The Robin coefficients \(\alpha_n\) used to define the vertex condition~\eqref{KirchhoffRobin}. To implement a Dirichlet vertex condition~\eqref{Dirichlet}, use a not-a-number (\NaN). Default value $0$ for Neumann-Kirchhoff conditions.
  \itemtt{y} The values of \(\psi\) at the vertices, set to \NaN\ on initialization. Used only for plotting.
\end{description}

We then examine the edges table, which has seven required fields:
\begin{lstlisting}[numbers=none]{matlab}
>> disp(G.Edges)
    EndNodes Weight   L    nx        x              y      Field7  
    ________ ______  ____  __  _____________  _____________  ____
     1  2       1      4    4  { 6x1 double}  { 6x1 double}    1
     2  2       1    6.28   8  {10x1 double}  {10x1 double}  0.79
\end{lstlisting}
The fields are:
\begin{description}
  \itemtt{EndNodes} This \(\abs*{\EE} \times 2\) array the initial and final vertices of the edges, i.e., the content of the input variables \source\ and \target.
  \itemtt{Weight} The weight \(w_j\) in vertex condition~\eqref{KirchhoffRobin}. Defaults to one if unset. 
  \itemtt{L} The array of edge lengths.
  \itemtt{nx} The number of discretization points on the interior of each edge. 
  \itemtt{x} The \((\mathtt{nx}+2)\) discretization points on each edge, including ghost points for the uniform discretization and vertices in the Chebyshev discretization.
  \itemtt{y} The value of \(\psi\) at the discretization points, initially set to \NaN.
  \itemtt{Field7} This field, named \listtt{integrationWeight}, is the spatial discretization \listtt{dx} for uniform discretizations and the Curtis-Clenshaw weights used to compute integrals with the Chebyshev discretization. 
\end{description}

Plotting coordinates are defined in a separate program, e.g., Fig.~\ref{fig:lasso}(a) was crearted using the program \listtt{lassoPlotCoords.m}. The coordinates are assigned to the graph \ttG\ with the command \listtt{G.addPlotCoords(@lassoPlotCoords)}
\noindent after \ttG\ has been created. Alternately, it can be included in the constructor by setting the optional argument \listtt{PlotCoordinateFcn} to the value \listtt{@lassoPlotCoords}.

 The plot coordinates are stored in fields \listtt{x1} and \listtt{x2} in both the \listtt{Edges} and \listtt{Nodes} tables. The MATLAB \listtt{plot} command is overloaded so  \listtt{G.plot}  plots the \listtt{y} coordinate over a skeleton of the graph in the \listtt{x1} and \listtt{x2} coordinates. Some graphs, such as those formed from the edges and vertices of a platonic solid, are best depicted in three space variables, so the user may define a third plot coordinate \listtt{x3}. If \listtt{x3} exists, the graph is plotted in three dimensions with the \listtt{y} coordinate represented by a color scale.

QGLAB includes templates for various commonly studied graphs and a template syntax that allows the quick creation of such graphs, such as the lasso and tetrahedron templates used in Sec.~\ref{sec:functions-and-plotting}. These have default parameters that can be overridden. A gallery of graph templates is included in the documentation.

\subsection{Basic operations}
\label{sec:operations}
A MATLAB live script is a rich document that includes runnable code and formatted text, entered with a simple word processor-like interface, and which integrates outputs including text and graphics, which can be exported to formats including PDF, \LaTeX, and HTML. QGLAB includes many examples created as live scripts and exported to HTML. Basic operations are described in the file \listtt{documentation/quantumGraphRoutines.mlx}.

\subsubsection{Function evaluation/plotting}
\label{sec:functions-and-plotting}
The command \listtt{applyFunctionToEdge} evaluates a function specified by a function handle, anonymous function, or constant value and assigns its value to edge \(\edge_j\). The command \listtt{applyFunctionsToAllEdges} applies a cell array of functions to all the edges; for example, the following commands define and plot a dumbbell quantum graph with the default parameters, plotted in Fig.~\ref{fig:firstExample}(a):

\vspace{.1cm}
\begin{lstlisting}[numbers=none]{matlab}
G=quantumGraphFromTemplate('dumbbell');
G.applyFunctionsToAllEdges({@sin,@(x)exp(-(x-2).^2),0});
G.plot
\end{lstlisting}
\begin{figure}[htbp] 
   \centering
   \begin{subfigure}[b]{0.4\textwidth}
\includegraphics[width=\textwidth]{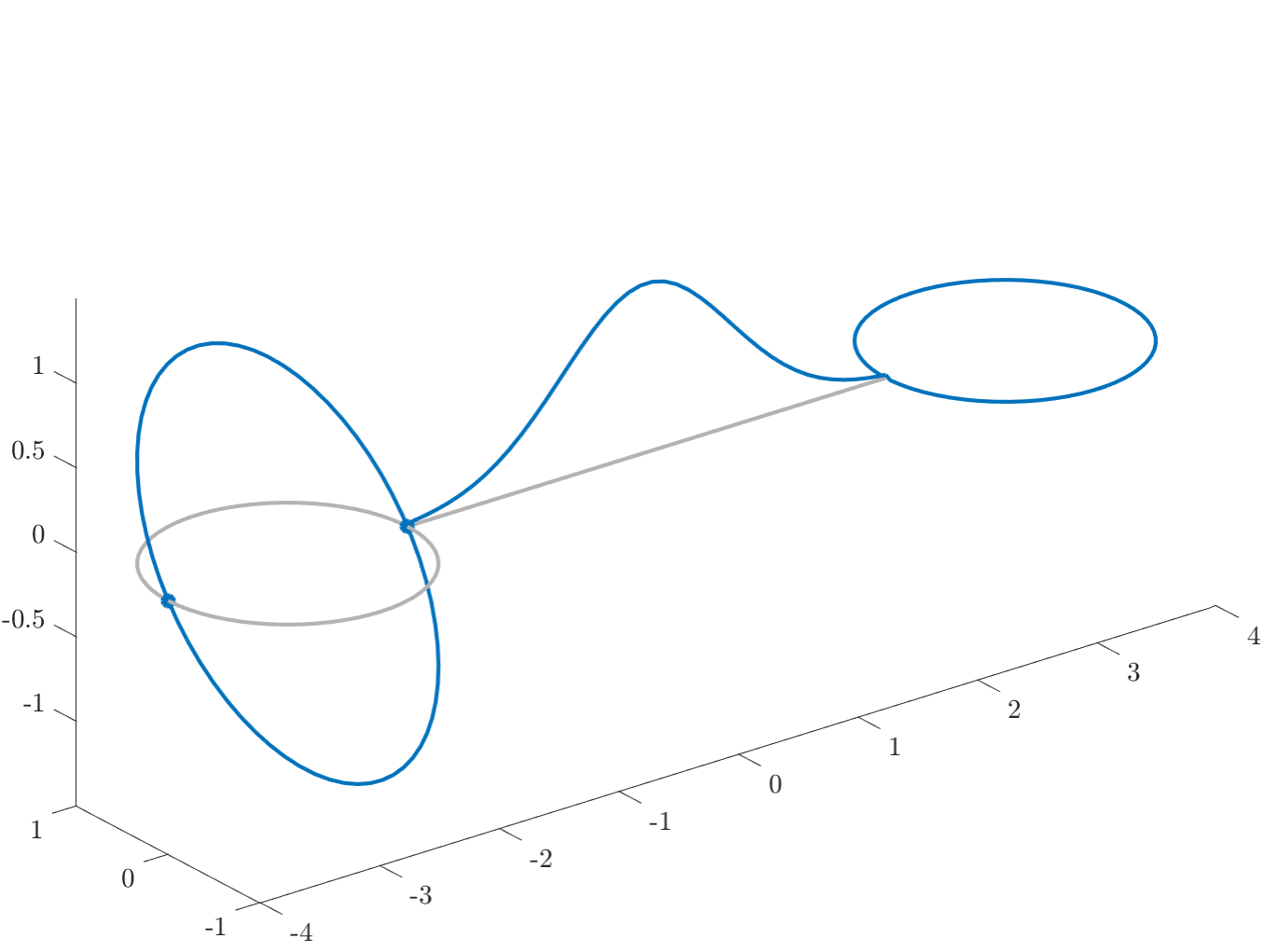}  
\caption{}
   \end{subfigure}
\begin{subfigure}[b]{0.4\textwidth}
    \includegraphics[width=\textwidth]{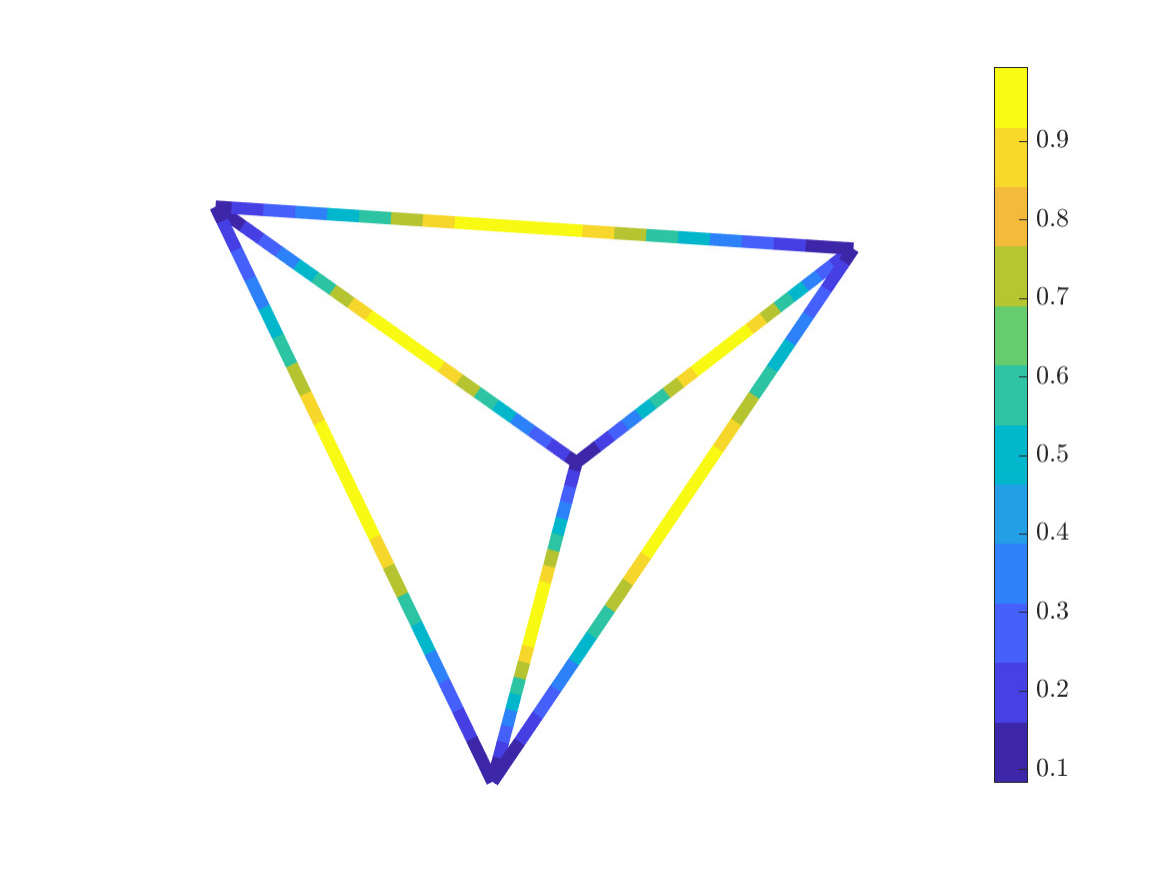}
    \caption{}
\end{subfigure}
   \caption{(a) A function defined on the edges of a dumbbell graph. (b) A function defined on the edges of a tetrahedral graph.}
   \label{fig:firstExample}
\end{figure}

QGLAB also provides some three-dimensional templates, for which the \listtt{y} values are plotted using a color scale, as in Fig~\ref{fig:firstExample}(b), where we plot Gaussians on all edges of a regular tetrahedron using the commands:
\begin{lstlisting}[numbers=none]{matlab}
G=solidTemplates('tetrahedron');
f=@(x)(exp((-10*(x-.5).^2))); 
G.applyFunctionsToAllEdges({f,f,f,f,f,f});
G.plot;
\end{lstlisting}

Finally, on some graphs with many edges, plots in three dimensions become a confusing tangle of curves, and it is more illuminating to visualize them as a color scale, using an overloaded \listtt{pcolor} command, demonstrated in Sec.~\ref{sec:function_listing}.

\subsubsection{Getting data on and off the graph}

The discretized numerical problems, including Poisson problem~\eqref{discretizedLaplacian} and the eigenvalue problem~\eqref{numericEvalProb}, are posed in terms of unknown column vectors. By contrast, this data is stored edge by edge in the \listtt{quantumGraph} object in \listtt{G.qg.Edges.y}. 
The command \listtt{graph2column} creates such a column vector from the data in graph \ttG, while the command \listtt{column2graph} loads the data from a column vector onto the edges of the graph including the vertices (Chebyshev) or ghost points (Uniform). Under the uniform discretization, the command also interpolates the data to the vertices. The \listtt{column2graph} command is also called by the command \listtt{G.plot(y)} to plot the contents of the vector \tty\ over the graph's skeleton.

\subsubsection{Other overloaded functions}

Many methods from MATLAB's directed-graph class have been overloaded so that, for example, a call to \listtt{G.numnodes} returns the number of nodes and \listtt{G.numedges} the number of edges. An overloaded \listtt{spy} command (along with additional formatting) was used to visualize the Laplacian matrix in Figs.~\ref{fig:lasso} and Fig.~\ref{fig:lassocheb}. Overloaded versions of the \listtt{norm} and \listtt{dot} commands are used frequently throughout the package. Sec.~\ref{sec:function_listing} of the \suppendix\ lists these functions.

\section{Extended examples}
\label{sec:examples}
Here, we give some brief examples of problem setup and accuracy. The \suppendix contains some additional details and other examples. A MATLAB live script for each, featuring plots and a convergence study, is included in the QGLAB GitHub repository and indexed in the {\bf readme} file~\cite{Goodman:2024}.

\subsection{The Poisson problem}
\label{sec:poisson_example}
We construct a discretization~\eqref{discretizedLaplacian} of the Poisson problem~\eqref{Poisson} using the graph in Fig.~\ref{fig:sillyGraph} with edges of lengths \( (\pi, 2\pi, 1, 2\pi, 2) \). The vertex conditions at \(\vertex_1\) and \(\vertex_2\) are Kirchhoff-Robin with \( \alpha_1= \alpha_2 = 1 \) and at \(\vertex_3\) the vertex condition is Dirichlet. The weight vector is $\vb{w}=(1,1,2,1,1)$. The potential vanishes except on edge \(\edge_1\) where $V_1=2\cos{2x}$. The nonhomogeneous terms are 
\[ f = \left(-\sin{3x},2\cos{2x},-4;-\sin{x},\sech{x}-2\sech^3{x} \right) \qand
 \psi =  ( 8, 3, \sech{2}).\]
The exact solution is
\[ \psi = \left(\sin{x},\sin^2{x},3x-2x^2,1+\sin{x},\sech{x}\right). \]

A minimal code to solve this problem is
\begin{lstlisting}[numbers=left,firstnumber=1,numberstyle=\tiny]{matlab}
s=[1 1 1 2 2]; t=[1 1 2 2 3]; L=[pi 2*pi 1 2*pi 2];
rc = [1 1 nan]; w = [1 1 2 1 1]; V = {@(x)(2*cos(2*x)),0,0,0,0};
G=quantumGraph(s,t,L,'RobinCoeff',rc,'Weight',w,'Potential',V);
f = G.applyFunctionsToAllEdges({@(x)-sin(3*x);@(x)2*cos(2*x);...
    -4;@(x)-sin(x);@(x)sech(x)-2*sech(x).^3});
phi = [8;3;sech(2)];
psi = G.solvePoisson('edgeData',f,'nodeData',phi);
\end{lstlisting}
Lines 1 and 2 define the quantum graph, the parameters that define the vertex conditions, and the potential. The variables \listtt{s}, \listtt{t}, and \listtt{L} define the source vertices, target vertices, and lengths of the five edges. The first two elements of \listtt{rc} define the Robin coefficients at vertices \(\vertex_{1,2}\), and the third component is \listtt{nan} (not-a-number), indicating the Dirichlet condition. The weight vector and the potential are given by \listtt{w} and \listtt{V}. Line 3 constructs the quantum graph from this data. It used the default centered-difference discretization and the default value \listtt{nX=20}, so it uses $h=\frac{1}{20}$ or a slightly smaller value on each edge to achieve an integer number of subintervals per edge. Lines 4-6 define the nonhomogeneous terms, and line 7 solves the discretized problem.

The maximum pointwise error of this solution is $1.02\times 10^{-3}$. Doubling the resolution by inserting the arguments \listtt{'nX',40} at Line 3 reduces the error to $2.56\times 10^{-4}$, a factor of $4.01$. This is consistent with the expected second-order convergence. Adding the arguments \listtt{'Discretization','Chebyshev'} at Line 3 switches to Chebyshev discretization. With 16 points per edge, this gives a maximum error of $1.80\times10^{-7}$; with 32, the error is $1.56\times10^{-12}$. This is consistent with spectral convergence.

\subsection{Eigenproblems}
\label{sec:eigen_example}

\begin{figure}[htbp] 
   \centering
   \includegraphics[width=0.5\textwidth]{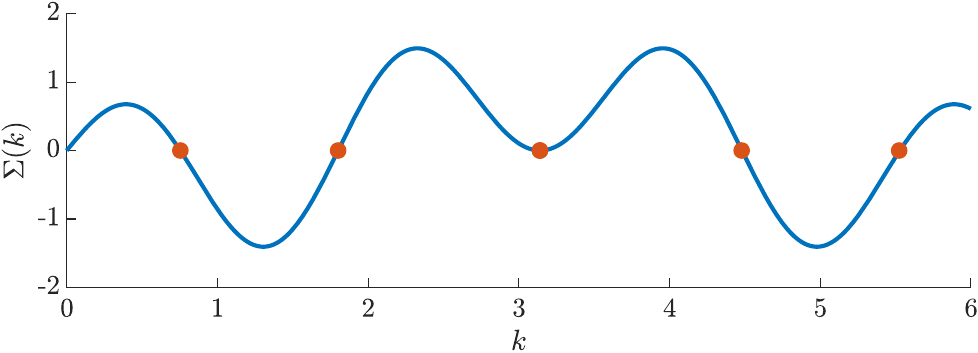} 
   \caption{The secular determinant \(\Sigma(k)\), of the Y-shaped graph discussed in the text, along with the computed values \(k_j=\sqrt{-\lambda_j}\), which sit right on top of the zeros. }
   \label{fig:Yshape_secdet}
\end{figure}

The example plotted in Fig.~\ref{fig:star_eigenfunctions} is computed in the live script \listtt{starEigenfunctionsDemo.mlx}. The Y-shaped graph has edges of lengths \(\left\{\frac{3}{2},1,1\right\}\). Its Kirchhoff conditions at both ends of the longer edge \(\edge_1\) and Dirichlet conditions at the remaining vertices are defined by setting the first two elements of the \listtt{robinCoeff} vector to zero the remaining two to \listtt{nan}:
\begin{lstlisting}[numbers=left,firstnumber=1,numberstyle=\tiny]{matlab}
G = quantumGraphFromTemplate('star','LVec',LVec=[1.5 1 1],...
                             'robinCoeff',[0 0 nan nan],'nX',40);
[V,lambda]=G.eigs(4); % Compute 4 eigenvalues
for k=1:4; G.plot(V(:,k));end;
\end{lstlisting}
Running \listtt{G.secularDet} returns the following secular determinant plotted in Fig.~\ref{fig:Yshape_secdet},
\[
\Sigma(k) = \frac{4}{3} \sin{\frac{k}{2}} (\cos{k}+1) \left(6 \cos^2{k}-3 \cos{k}-1\right).
\]

The $\Sigma(k)$ plot shows that the eigenvalue \(\lambda=-\pi^2\approx -0.9865\) has multiplicity two. Generally, a numerical eigenvalue solver will return two closely spaced eigenvalues rather than a double eigenvalue. The graph \(\Gamma\) is symmetric under the interchange of the edges \(\edge_2\) and \(\edge_3\). The multiplicity-one eigenfunctions respect this symmetry, but the multiplicity-two eigenfunctions returned by \eigs\ do not. The script takes appropriate linear combinations of the two computed eigenvectors to produce odd and even eigenvectors with respect to this symmetry; see panels 3 and 4 of Fig.~\ref{fig:star_eigenfunctions}.  Sec.~\ref{sec:eigen_suppendix} of the \suppendix\ discusses the convergence, which is standard.

\subsection{Time-dependent problems}
\label{sec:NLS_example}
The nonlinear Schrödinger equation~\eqref{NLS} is the prototype of a dispersive nonlinear PDE and many studies have considered its evolution on quantum graphs. Kairzhan et al.\ consider the cubic NLS equation posed on a "star graph" consisting of three half-lines joined at a single vertex~\cite{Kairzhan:2019fe}. The evolution conserves mass, i.e., the squared \(L^2\) norm~\eqref{Lpnorm}, and the energy~\eqref{energyNLS}, but in general, does not conserve momentum. If, however, the parameters \(w_m\) are appropriately chosen in the weighted Kirchhoff-Robin vertex condition~\eqref{KirchhoffRobin} to form a so-called "balanced" star graph, and the initial condition is chosen to lie in a particular invariant subspace, then the dynamics on this graph do conserve momentum.

We simulate the collision of a soliton propagating on edge \(\e_1\) toward the vertex results in Fig.~\ref{fig:NLSonStar} using the implicit/explicit Runge-Kutta solver \listtt{qgdeSDIRK443}. The first simulation is computed on a  "balanced" graph whose vertex condition is defined by the weight vector \(\vb{w} = (2,1,1)\) in the first line of the following script.
\begin{lstlisting}[numbers=left,firstnumber=1,numberstyle=\tiny]{matlab}
G =quantumGraphFromTemplate('star','LVec',30,'weight',[2 1 1]);
soliton =@(x,t,v,x0)exp(1i*(-v*x/2-(1-v^2/4)*t)).*sech(x-x0-v*t);
init1 = @(x)soliton(x,0,-2,15); init2 =@(x)soliton(x,0,2,-15);
u0 = G.applyFunctionsToAllEdges({init1,init2,init2});
mu = -1i; F =@(z) -2i * z.^2.* conj(z);
tFinal = 11; dt = 0.01; tPrint =0.5; nSkip = tPrint/dt;
[t,u] = G.qgdeSDIRK443(mu,F,tFinal,u0,dt,'nSkip',nSkip);
\end{lstlisting}
The soliton splits into two, each new soliton propagating along the edge with its original amplitude and velocity. At this discretization, the mass, energy, and momentum are all conserved to 5, 4, and 4 digits, respectively.

We next set \(\vb{w} = (1,1,1)\). Now, a significant fraction of the soliton is reflected and propagates backward along the incoming edge. Mass and energy conservation are slightly worse, but momentum changes by $O(1)$. The standard test of computing the numerical solution with time steps $\tau$, $\tau/2$, and $\tau/4$ shows that the difference between subsequent solutions decreases about eight times, indicating third-order convergence.

\section{Conclusions}
\label{sec:conc}

QGLAB is a robust and versatile MATLAB package for computing solutions to the spectral accuracy of linear and nonlinear problems on Quantum Graphs. It allows users to build graph models quickly, analyze their spectrum, compute nonlinear bifurcations, and solve evolution equations. The algorithms are implemented at a high level, hiding most implementation details and allowing the user to focus on the mathematical problem, not the numerical and algorithmic details.

Linear and nonlinear PDEs on quantum graphs remain a vibrant area of analysis in spectral geometry in which the interaction of geometry, topology, and symmetry gives rise to diverse mathematical questions~\cite{alon2018nodal,berkolaiko2019surgery,gnutzmann2006quantum} with many open problems left to explore. Many previously studied problems on combinatorial graphs have analogies on metric graphs that remain open and where the spectrum of behaviors is likely to be much richer. For example, the spectral optimization of combinatorial graphs has been studied in~\cite{osting2017spectrally}, and others have examined how the symmetries of discrete Laplacians can lead to interesting spectral features such as Dirac points and flat bands~\cite{lim2020dirac,morales2016simple}. The study of time-dependent evolution equations on quantum graphs remains in its infancy~\cite{Dutykh:2018,Kairzhan:2019fe}. 
QGLAB is an ideal tool for exploring these problems. 

\newpage
\appendix

\begin{center}
    {\bf Supplementary Material}
\end{center}

This \suppendix\ contains two sections. The first, Section~\ref{sec:appendix_examples}, is devoted to demonstrating both the implementation and efficacy of QGLAB on a variety of examples, including stationary problems---eigenvalue problems, the Poisson equation, and the computation and continuation of standing waves---in Section~\ref{sec:stationary_examples} and evolutionary PDE problems in Section~\ref{sec:evolution_examples}. All the examples are included as live scripts (MATLAB \listtt{.mlx} files) in the directory \listtt{source/examples}. The second part, Sec.~\ref{sec:function_listing}, contains a complete listing of user-callable function definitions and explicit instructions for their use. 

\section{Extended examples}
\label{sec:appendix_examples}

\subsection{Stationary problems}
\label{sec:stationary_examples}

\subsubsection{Eigenproblems}
\label{sec:eigen_suppendix}
Here, we report in greater detail on the accuracy of the eigenproblem calculated in Sec.~\ref{sec:eigen_example}. We find the exact eigenvalues from the zeros of the secular determinant. Then we compute the finite-difference approximation with \(h=\frac{1}{40}\) and \(h=\frac{1}{80}\). The ratio of these is about 4, which shows the method is second order. Finally, we compute the same eigenvalues using the Chebyshev discretization with 30 points on the long edge and 20 on the short edges. Since all errors are less than \(10^{-10}\), we conclude the accuracy is spectral.

\begin{center}
\begin{tabular}{rrcccc}
$k=\sqrt{-\lambda}$ & $\lambda$ & err$_{\tau=\frac{1}{40}}$  & err$_{\tau=\frac{1}{80}}$ & ratio & err$_{\rm Cheb}$\\
\hline \\
$\cos^{-1}\left(\frac{1}{12} \left(\sqrt{33}+3\right)\right)$ & -0.569 & 1.687e-05 & 4.216e-06 & 4.000 & 2.195e-12\\ 
$\cos^{-1}\left(\frac{1}{12} \left(\sqrt{33}-3\right)\right)$ & -3.246 & 5.486e-04 & 1.372e-04 & 4.000 & 1.177e-12\\ 
$\pi$& -9.870 & 5.072e-03 & 1.268e-03 & 3.999 & 1.506e-11\\ 
$\pi$ & -9.870 & 5.072e-03 & 1.268e-03 & 3.999 & 3.020e-14\\ 
$2\pi-\cos^{-1}\left(\frac{1}{12} \left(\sqrt{33}-3\right)\right)$ & -20.085 & 2.100e-02 & 5.252e-03 & 3.999 & 3.830e-11\\ 
\(2\pi-\cos^{-1}\left(\frac{1}{12} \left(\sqrt{33}+3\right)\right)\) & -30.568 & 4.864e-02 & 1.216e-02 & 3.998 &1.669e-11 \\ 
\(2\pi\) & -39.4784 & 8.111e-02 & 2.029e-02 & 3.998 & 1.847e-13
    \end{tabular}
\end{center}

\subsubsection{Nonlinear standing waves and bifurcation diagrams}
\label{sec:nonlinear_standing_examples}

\paragraph{\textbf{Computing individual solutions}}
We begin with an example computing a single solution to the stationary cubic NLS~\eqref{stationary_NLS} on a dumbbell graph:
\begin{lstlisting}[numbers=left,firstnumber=1,numberstyle=\tiny]{matlab}
G = quantumGraphFromTemplate('dumbbell');
fcns = getNLSFunctionsGraph(G);
Lambda = -1;
f = @(z)fcns.f(z,Lambda); M = @(z)fcns.fLinMatrix(z,Lambda);
y0 = G.applyFunctionsToAllEdges({0,@(x)sech((x-2)),0});
y = solveNewton(y0,f,M); G.plot(y)
\end{lstlisting}
The function \listtt{getNLSFunctionsGraph} defines the discretized version of the nonlinear functional and several of its partial derivatives and assigns them to a structure array called \listtt{fcns}. By default, this uses the function \(f(z) = 2 z^3\) from Eq.~\eqref{stationary_NLS}. The user may provide a symbolic function of one variable as an optional argument, and MATLAB will compute all the required partial derivatives symbolically. The Newton-Raphson solver that is iterated to solve the system requires both the functional and its linearization with respect to \(\Psi\). These are stored in two fields \listtt{fcns.f} and \listtt{fcns.fLinMatrix}, which are functions of two inputs \listtt{z} and \listtt{Lambda}. The continuation algorithm considers Eq.~\eqref{stationary_NLS} as a function of both \(\Psi\) and \(\Lambda\), but in this first example, we fix \(\Lambda=-1\) and consider only \(\Psi\) as unknown. In line 4, \emph{anonymous functions} are used to instruct MATLAB to consider them as functions of \(\Psi\) alone. We search for a unimodal solution to Eq.~\eqref{stationary_NLS} with \(\Lambda=1\) centered on the central edge of a dumbbell graph, so we prepare an initial guess in line 5 consisting of a hyperbolic secant centered on the central edge and zeros on the two looping edges. The \listtt{solveNewton} command finds the standing wave. The result of the \listtt{plot} command is shown in Fig.~\ref{fig:stationaryNLS}(a).

For graphs with a large number of edges, generating an initial guess with the approach of line~6 would be impractical, so QGLAB provides a convenient function \listtt{applyGraphicalFunction} which applies a function to the coordinate functions used to plot the graph. In Fig.~\ref{fig:stationaryNLS}(b), we find a standing wave on a spiderweb graph, found in the QGLAB template library, using as an initial guess the function \(\sech{(r)}\) where \(r\) is the Euclidean distance from the central point to a point on the graph as laid out in two dimensions.

\paragraph{\textbf{Continuation of solutions}}
We can learn more about the stationary problem by considering branches of standing waves and their bifurcations than by computing individual solutions. Well-established and sophisticated software packages for such computations include AUTO and MatCont for ODE systems and pde2path for elliptic PDE~\cite{Dhooge:2003vy, Dhooge:2008gq, Doedel:2007, Uecker:2014}. The capabilities of QGLAB are much more modest but allow for the simple setup and solution to continuation and bifurcation problems on quantum graphs, following branches around folds, detection of bifurcation points, and changing branches at such points. 

\begin{figure}[htbp] 
   \centering
      \begin{subfigure}[b]{0.4\textwidth}
          \includegraphics[width=\textwidth]{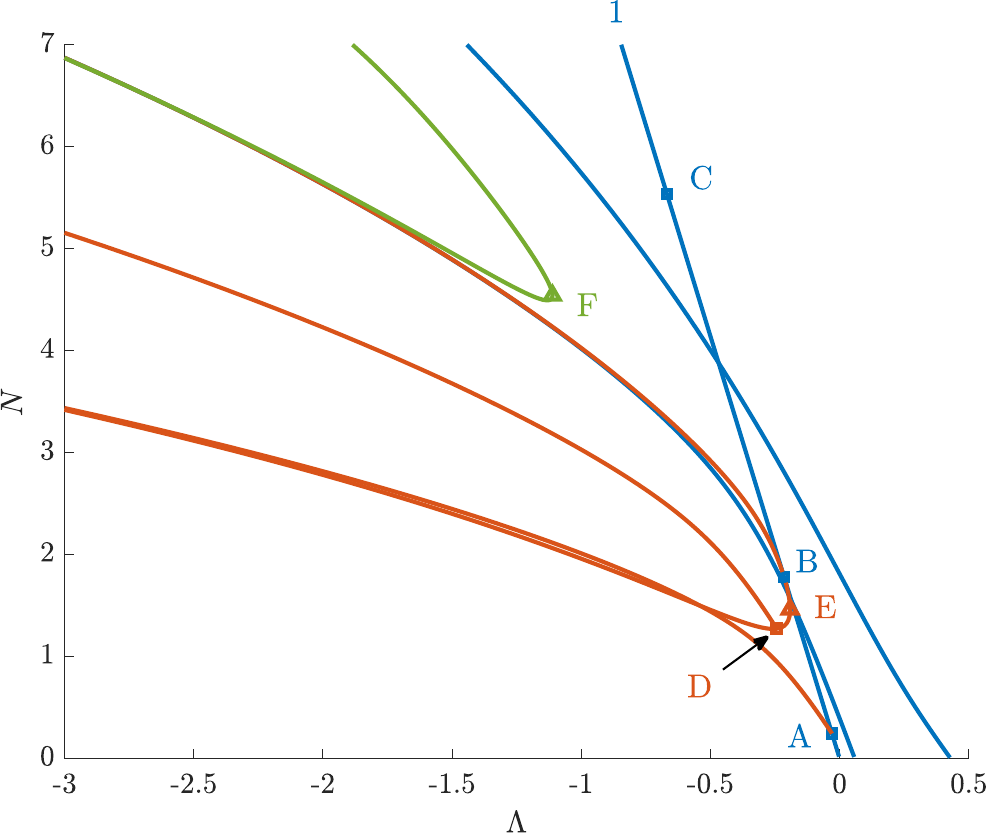}
          \caption{}
      \end{subfigure}
\begin{subfigure}[b]{0.4\textwidth}
     \includegraphics[width=\textwidth]{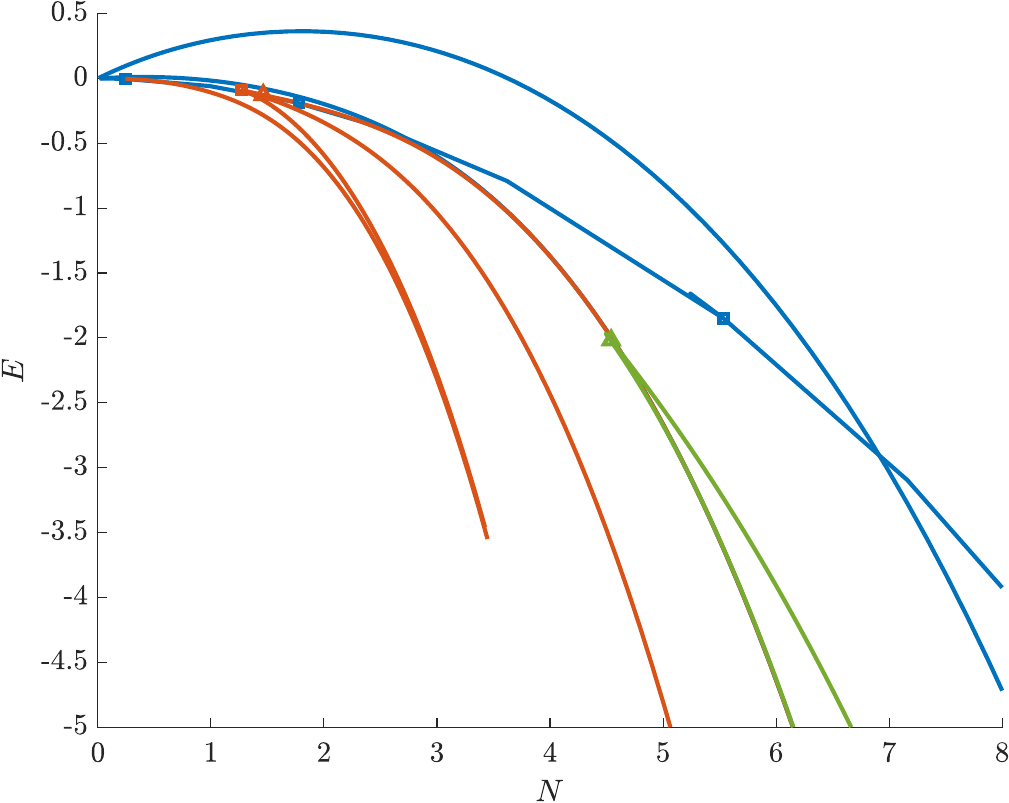}
     \caption{}
\end{subfigure}
   \caption{(a) A partial bifurcation diagram for the dumbbell graph. The three blue curves are the continuations of linear eigenfunctions. The red curves were computed by continuing from branching bifurcations. The green curve was computed by computing a single large amplitude solution and then continuing it. Branching bifurcations marked with squares and folds with triangles. (b) The same diagram, plotted in different variables.}
   \label{fig:dumbbellBifurcation}
\end{figure}

An extended example of numerical continuation is presented in the live script that is titled \listtt{continuationInstructions.mlx}, which presents a computation of a partial bifurcation diagram of the cubic NLS equation on a dumbbell graph in Fig.~\ref{fig:dumbbellBifurcation}, reproducing a figure from~\cite{Goodman:2019go}, which contains far more details and graphs of several of the solutions at various points on the bifurcation diagram. 

This figure comprises nine separately-computed curves, each representing dozens of solutions to Eq.~\eqref{stationary_NLS}. The curves were initialized in three different ways. The first type, plotted in blue, consists of nonlinear continuations of linear eigenfunctions. 
We have plotted three such branches but focus on the branch labeled \textbf{1}. 
This branch represents the nonlinear continuation of the null eigenvector of the Laplacian on this quantum graph.  The value of \(\Psi\) is constant on all solutions on this branch, with
\begin{equation} \label{Psiconstant}
\Psi= \sqrt{\frac{-\Lambda}{2}}. 
\end{equation}

It is straightforward to show that if \(\lambda\) is an eigenvalue of the operator \(-\Laplacian\), then branch \textbf{1} has a bifurcation point at \(\Lambda = -\lambda/2\)~\cite{Goodman:2019go,Marzuola:2016bl}. 
QGLAB automatically computes the direction in which branches fork from bifurcation points, and the diagram shows two families that emerge from such points. At the points marked \textbf{A}, \textbf{B}, and \textbf{C}, QGLAB has detected bifurcation points on branch \textbf{1}, and we have chosen to follow the first two.
The branch that bifurcates from branch \textbf{A}, which seems to intersect branch \textbf{1} transversely, is a pitchfork bifurcation, while the branch that bifurcates from \textbf{B} tangentially to branch \textbf{1} and extends in both directions is a transcritical bifurcation. This last branch itself has a limit (fold) point at \textbf{E} and a pitchfork bifurcation at \textbf{D}. The final branch, plotted in green, was generated by first computing a single high-frequency bifurcation with large amplitude pulses on the dumbbell handle and one ring, saving it to a file, and then continuing that solution.

QGLAB stores all the data for branches, bifurcation points, and individual solutions logically and hierarchically and has routines for retrieving and plotting individual solutions and curves of solutions so that the user can largely avoid low-level interactions with the data. By default, it plots the frequency of standing waves versus their power, but it can also plot the energy~\eqref{energyNLS}, as shown in the right image of Fig.~\ref{fig:dumbbellBifurcation}.

\begin{figure}
    \centering
    \includegraphics[width=0.6\textwidth]{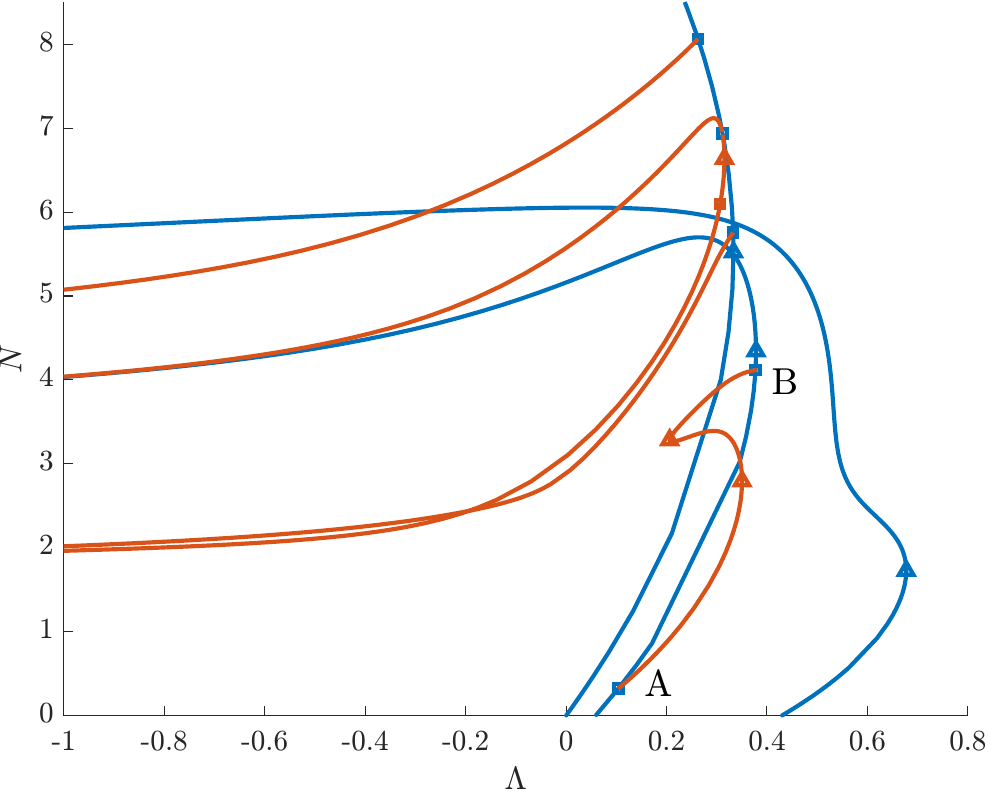}
    \caption{A partial bifurcation diagram of the stationary NLS equation on a dumbbell quantum graph with a cubic-quintic nonlinearity.}
    \label{fig:bifurcation35}
\end{figure}

The nonlinear term in stationary NLS~\eqref{stationary_NLS} can be changed by simply changing the definition of \(f(z)\) to any analytic function satisfying \(f(0)=0\) (so that the linearization at zero remains unchanged and the continuation of linear eigenfunctions from zero can be easily computed). In the example \listtt{dumbbellcontinuation35.mlx}, we change the right hand side to \(f(z)=-2z^3 + 3z^5\) which is defocusing for small values of \(\abs{z}\) and focusing for large values. A partial bifurcation diagram for this system is shown in Fig.~\ref{fig:bifurcation35}, consisting of three branches that bifurcate from zero in the direction of the eigenfunctions, albeit with a frequency that initially increases with increasing power before changing direction and decreasing. The leftmost branch remains constant in space, and its power increases monotonically along the branch. In contrast, the other two branches have decreasing power as the frequency decreases past a certain point. 

Especially interesting is the branch that bifurcates from the point  \textbf{A} on the middle branch. This middle branch is the continuation of the first excited eigenfunction, which has an odd symmetry about the central point on the dumbbell. At this point, we find a symmetry-breaking pitchfork bifurcation, with two asymmetric branches related by a reflection symmetry. This asymmetric branch continues to the point \textbf{B}, at which point it collides again with the same branch from which it bifurcated at \textbf{A} and begins retracing its original path. This branch traces out a closed curve in solution space, with the sign of the perturbation term flipping each time the branch passes the bifurcation points. Thus, we instructed the continuation program to stop after a finite number of points on the curve are computed by setting the parameter \listtt{maxPoints} as described in Sec.~\ref{sec:appendixContinuation}.

An advantage of the continuation/bifurcation approach is that it illuminates how branches relate to each other. This is well illustrated using the example of a "necklace" quantum graph, also considered by Besse et al.~\cite{Besse:2021}. This graph consists of loops alternating with single edges. The necklace graph shown above in Fig.~\ref{fig:necklace}(a) consists of 54 such alternating pairs, with segments of length \(1\) and pearls comprised of two edges, each of length \(\pi/2\). Fig.~\ref{fig:necklace}(b) shows a partial bifurcation diagram for the focusing cubic NLS equation on this graph.

We focus on branch \textbf{1} and a few branches arising from bifurcations from this branch and its descendants. As in the first example, the constant-valued solution on this branch satisfies Eq.~\eqref{Psiconstant}, and bifurcations occur where the frequency is half of an eigenvalue of the linear problem. However, this eigenvalue has a geometric multiplicity of two in this case. In bifurcation theory, the system is said to undergo a \emph{codimension-two} bifurcation at this point. QGLAB has not implemented methods for detecting higher codimension bifurcation points and calculating branches emanating from bifurcations of codimension two or higher. Such methods exist and are implemented in the packages cited above; an approach that obviates the need to calculate higher-order normal forms is the deflated continuation method due to Farrell and collaborators~\cite{Farrell:2015}. 

\begin{figure}[htbp] 
   \centering
   \includegraphics[width=0.8\textwidth]{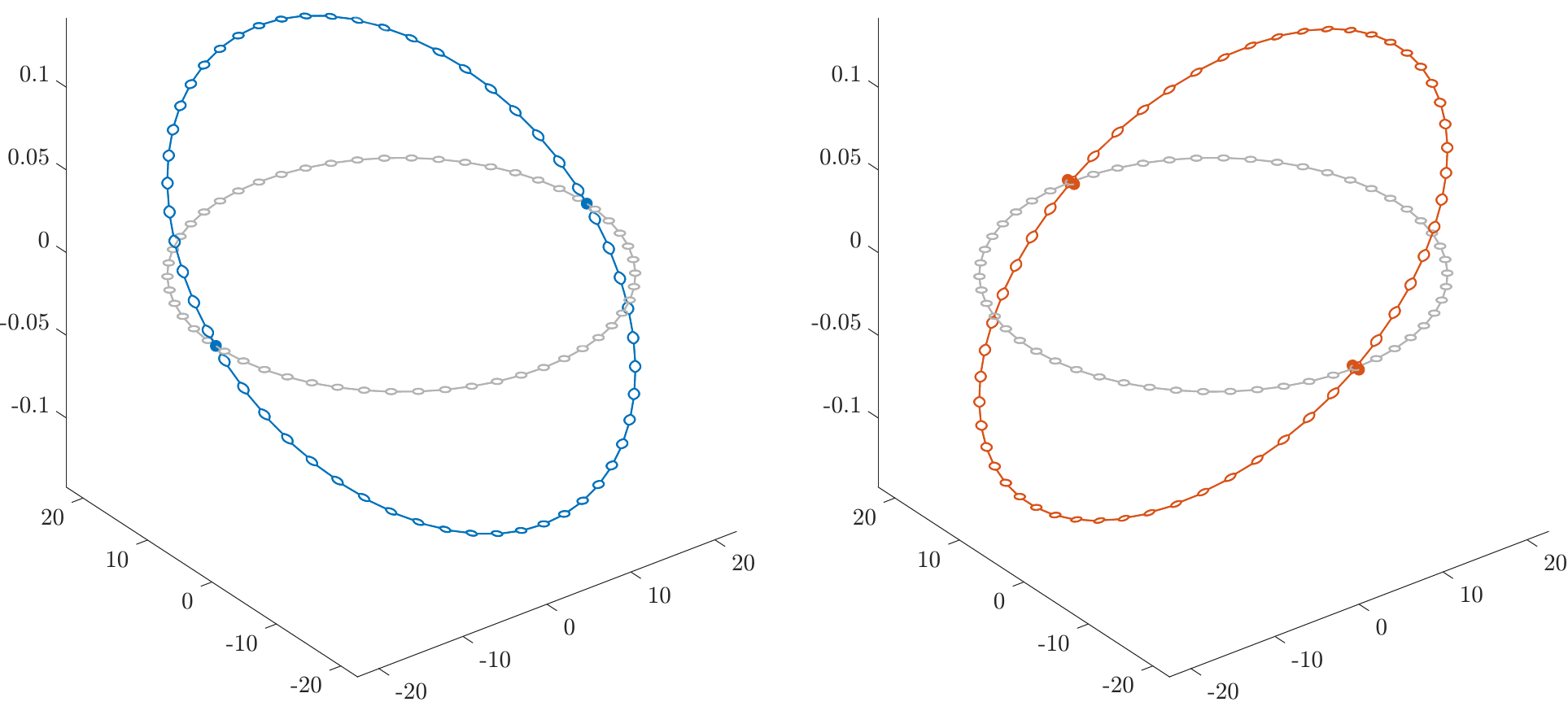} 
   \caption{The eigenfunctions corresponding to the smallest nonzero eigenvalue on the necklace quantum described in the text. The left eigenfunction has two nodes on ``strings'' and four local extrema on ``pearls'', while the right eigenfunction has four nodes on ``pearls'' and two local extrema on ``strings.''}
   \label{fig:necklace_eigenfunctions}
\end{figure}

The double-zero eigenvalue at this bifurcation has two orthogonal eigenfunctions plotted in Fig.~\ref{fig:necklace_eigenfunctions}. These may be thought of as the analog of the sine and cosine modes of the second derivative operator on the circle. While any linear combination of these two eigenfunctions is also an eigenfunction, we have chosen the two modes so that one has its maximum at the center of a single strand and the other at the center of a double strand. The nonlinear standing waves that bifurcate from branch \textbf{1} at the point \textbf{A} do so in the direction of these two eigenfunctions. Close to the bifurcation, the two solution curves are indistinguishable when plotted in these coordinates but separate for more negative frequencies. The standard algorithm that QGLAB uses to detect bifurcations works not by computing all the eigenvalues of the linearization and counting their eigenvalues, which would be slow, but by efficiently calculating the sign of the associated determinant using an \emph{LU}-decomposition and detecting when it changes. This works efficiently at codimension-one bifurcations but fails at codimension-two bifurcations like this one. As this would predict, the algorithm that detects bifurcations fails to find a bifurcation at \textbf{A} and does not compute the branching direction.

The branches~\textbf{2} and~\textbf{3} are calculated by first computing a single standing wave with frequency \(\Lambda=-4\) and either a single \(\sech\)-like hump centered on a string or two \(\sech\)-like humps centered on the two edges on the pearl and then continuing the branches toward the bifurcation point~\textbf{A}. Branch~\textbf{4} bifurcates from branch~\textbf{3} at the point~\textbf{B}, breaking the symmetry between the two edges of the pearl. By plotting this bifurcation diagram in the same coordinates as in the right image of Fig.~\ref{fig:dumbbellBifurcation}, we confirm the statement of Ref.~\cite{Besse:2021} that this branch represents the ground state at large amplitude. At point~\textbf{C}, Branch~\textbf{3} undergoes a second symmetry-breaking bifurcation, giving rise to branch~\textbf{5}, on which the two-humped standing wave on the pearl moves from the center of the pearl's edges toward either vertex. A similar symmetry bifurcation occurs on Branch~\textbf{2} at point~\textbf{D}, giving rise to Branch~\textbf{6}, along which the standing wave on the string moves away from the string's center and toward a vertex. Branches~\textbf{5} and~\textbf{6} appear to converge as \(\Lambda\) is further decreased. Representative standing waves along these five branches of the bifurcation diagram at \(\Lambda \approx -4\) are shown in Fig.~\ref{fig:necklace}(c) above.

Finally, conducting a proper continuation study of standing waves on an infinite necklace is difficult. For a fixed number of pearls, the total width of the standing wave is restricted by the circumference, but in the infinite limit, branches~\textbf{2} and~\textbf{3} bifurcate not from the solution of constant amplitude, but from the zero solution, with a width that diverges as the amplitude goes to zero. The limiting behavior exists for the standing waves of the standard cubic NLS problem. However, in that case, a standard method allows the width of the interval to increase, namely using a non-uniform discretization that widens to accommodate the slowing spatial decay rate. Such a trick is unavailable on the quantum graph, where the length scale imposed by the graph's edges precludes this approach.

\subsection{Evolutionary PDE}
\label{sec:evolution_examples}

We now discuss the full MATLAB implementation of the heat and sine-Gordon examples constructed in the article.

\subsubsection{The heat equation}
\label{sec:heat_example}
\begin{figure}[htbp] 
   \centering
   \includegraphics[width=0.65\textwidth]{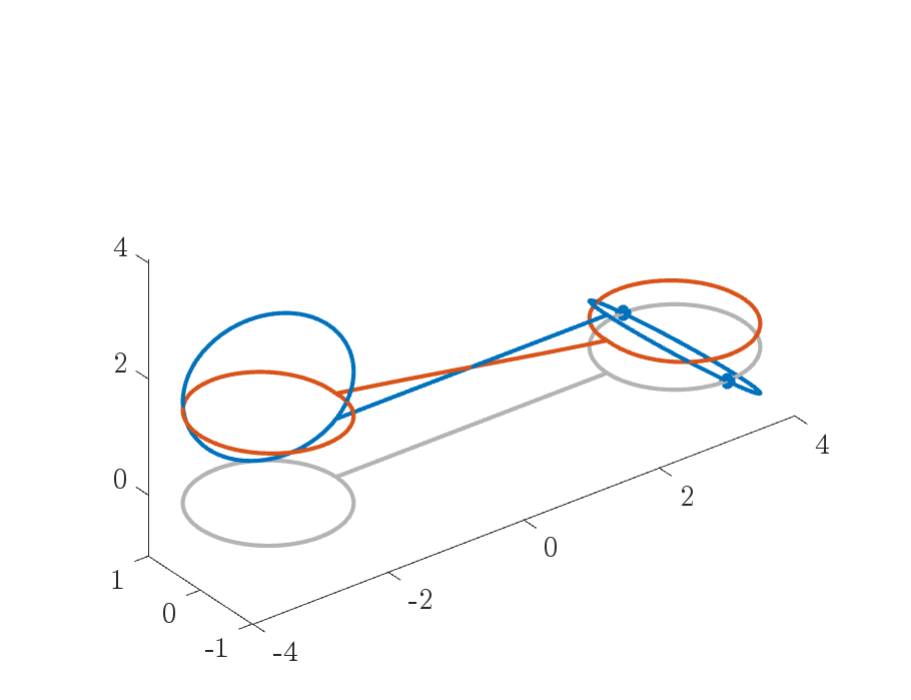} 
   \caption{The initial (blue) and final (red) states of the heat equation on a dumbbell graph computed using the Crank-Nicholson code in the text.}
   \label{fig:heatOnDumbbell}
\end{figure}

In Sec.~\ref{sec:elementary_time} we derived Eq.~\eqref{CNstepper} to evolve the solution of the heat equation over one time step. We apply this code to a dumbbell graph in the live script \listtt{heatOnDumbbell}. After removing the code for plotting and calculating the conserved total heat, the code reads
\begin{lstlisting}[numbers=left,firstnumber=1,numberstyle=\tiny]{matlab}
G = quantumGraphFromTemplate('dumbbell');
y=G.applyFunctionsToAllEdges({@(x)(2-2*cos(x-pi/3)),1,@cos});
dt=0.01; tFinal=10; nStep=tFinal/dt;
L0 = Phi.laplacianMatrixWithZeros;
P0 = Phi.interpolationMatrixWithZeros;
LVC = Phi.laplacianMatrixWithVC;
PVC = Phi.interpolationMatrixWithVC;
LPlus =  P0  + (h/2)*L0;
LMinus = PVC - (h/2)*LVC;
for k=1:nStep
    y = LMinus \ (LPlus*y);
end
\end{lstlisting}
This solution's initial and final states are shown in Fig.~\ref{fig:heatOnDumbbell}. The total heat is conserved to twelve digits by this calculation.

\subsubsection{The sine-Gordon equation}
\label{sec:sG_example}
\begin{figure}[htbp] 
   \centering
   \includegraphics[width=\textwidth]{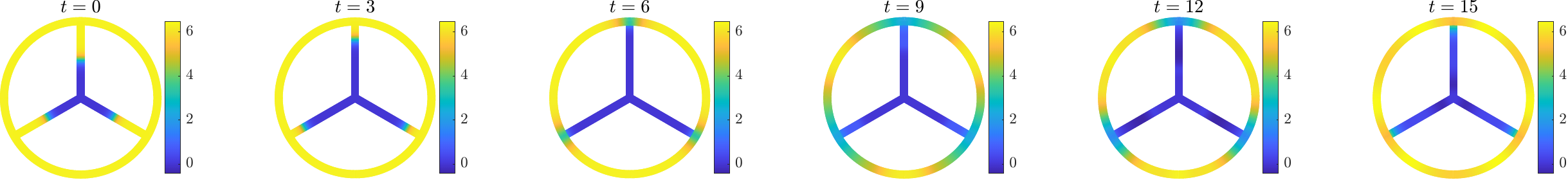}
   \includegraphics[width=\textwidth]{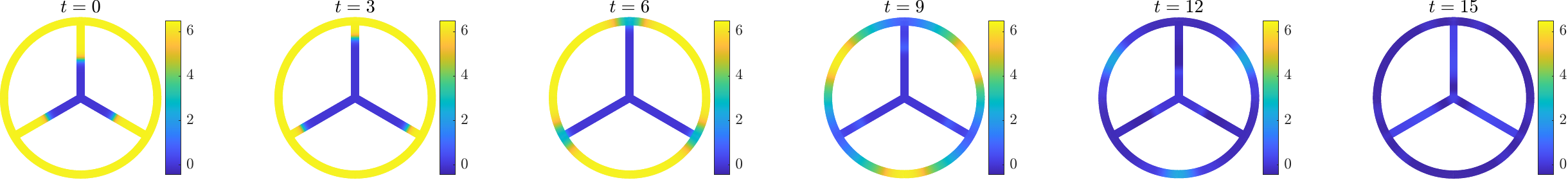}
   \caption{Evolution of sine-Gordon solitons propagating along the edges of a tetrahedron (deformed for plotting). (Top) the vertices reflect solitons with \(c=0.9\) while (Bottom) those with \(c=0.95\) are transmitted.}
   \label{fig:sGonTetra}
\end{figure}

The sine-Gordon equation on the line supports solitons, traveling solutions of the form
\[
\psi(x,t) = 4 \tan^{-1}\left( e^{(x-ct)/\sqrt{1-c^2}}\right), \qqtext{for any} -1<c<1.
\]
Following~\cite{Dutykh:2018}, we initialize kinks on three edges of the graph formed by the edges of a regular tetrahedron, heading away from their common vertex. We consider two initial conditions: the first with \(c=0.9\) and the second with \(c=0.95\). These are plotted in Fig.~\ref{fig:sGonTetra}, with the tetrahedron flattened into the shape of a wheel with three spokes (thus, distance in the plot does not uniformly represent distance on the metric graph). The top row shows the first case, in which the three solitons are reflected after encountering vertices, while in the second case, the faster solitons can pass through the vertices.

\section{Function Listing and Detailed Instructions}
\label{sec:function_listing}
\begin{figure}
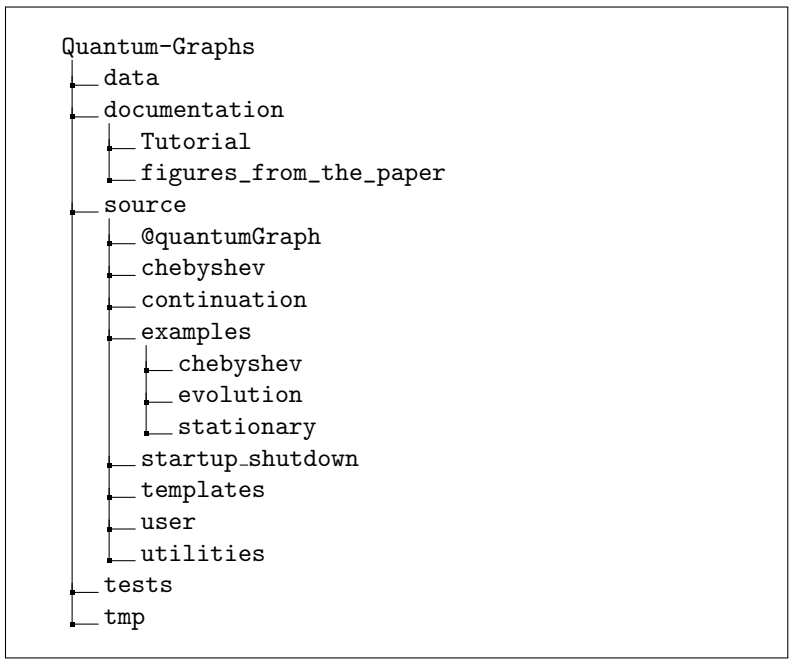

\begin{framed}[0.8\linewidth]
\begin{minipage}[t]{0.6\linewidth}
\dirtree{%
.1 Quantum-Graphs.
.2 data.
.2 documentation.
.3 Tutorial.
.3 figures_from_the_paper.
.2 source.
.3 @quantumGraph.
.3 chebyshev.
.3 continuation.
.3 examples.
.4 chebyshev.
.4 evolution.
.4 stationary.
.3 startup\_shutdown.
.3 templates.
.3 user.
.3 utilities.
.2 tests.
.2 tmp.
}
\end{minipage}
\end{framed}
\caption{Directory structure of QGLAB}
\label{fig:path}
\end{figure}

QGLAB is implemented as a MATLAB \emph{Project}. After starting MATLAB, the user should open the folder titled \listtt{Quantum-Graphs}, whose subfolder structure is shown in Fig.~\ref{fig:path}. 
Among the files listed in the MATLAB Desktop's \listtt{Current Folder}\ pane is the \emph{project file} \listtt{QGobject.prj}, which can be opened by double clicking. This opens the Project Window, adds the necessary QGLAB directories to MATLAB's search path, and changes the plotting preferences needed to render the graphics correctly. To end the QGLAB session, close the Project Window or quit MATLAB. This will remove the QGLAB directories from the search path and restore the user's default plotting preferences, which are held in the folder \listtt{tmp} while \listtt{QGLAB} is running.

The MATLAB code is contained in the subfolders of the folder \source. Most importantly, the folder \listtt{@quantumGraph} contains the \emph{constructor} file \listtt{quantumgraph.m}, which defines the class and initiates an instance, as well as all the class methods, i.e., the functions that act on quantum graph objects. As their first input argument, all MATLAB methods must have a \qg\ object \ttG. For example, the overloaded eigensolver method \eigs\ is defined as
\listtt{function [v,d]=eigs(G,n)}, where \listtt{n} is the number of eigenvalues to calculate. It can be called using either the standard function syntax \listtt{[v,d]=eigs(G,n)} or with the preferred syntax for methods \listtt{[v,d]=G.eigs(n)}. 

\subsection{The Quantum Graph Constructor}
\label{sec:constructor}
The first step to working with QGLAB is initializing a quantum graph object using its constructor function titled \qGraph. As detailed in Sec.~\ref{sec:initializing}, it takes three required arguments
\begin{description}
 \itemtt{source} and \target\ are two vectors of positive integers. The entries \listtt{source(m)} and \listtt{target(m)} represent the initial and final nodes of the edge $\edge_m$. Thus, these two vectors must be of the length $\abs{\EE}$ and each integer $m$ satisfying $1\le m \le \abs{\VV}$ must appear in at least one of the two vectors to guarantee that the graph is connected. MATLAB's \digraph\ constructor automatically sorts the edges to avoid confusion, \qGraph checks to make sure the edges are sorted the same way and throws an error if they are not.
 \itemtt{L} May be either a positive real number or a vector of length $\abs{\EE}$ of positive real numbers. If \ttL\ is scalar, the constructor assumes all edges are the same length.
 \end{description}
 It also may take the following optional arguments
 \begin{description}
 \itemtt{Discretization} One may take the values \Uniform\ (default), \Chebyshev, or \none. If \none, then no discretization is constructed, and the only available method, besides simple methods that query the graph's properties, is \listtt{secularDet}, which computes the secular determinant.
 \itemtt{nxVec} Defines the number of points used to discretize the edges. A vector value gives the number of discretization points on each edge, but if scalar, its behavior depends on the discretization; if \Uniform, then it gives the approximate number of points per unit edge length, while if \Chebyshev, then it gives the number of discretization points on each edge. Default: 20.
 \itemtt{RobinCoeff} The vector of Robin coefficients $\alpha_n$ in Eq.~\eqref{KirchhoffRobin}. Use the value \NaN\ to indicate the Dirichlet boundary condition~\eqref{Dirichlet}. If scalar, apply the same value at all vertices. Default: 0.
 \itemtt{Weight} The vector of weights $w_m$ in Eq.~\eqref{KirchhoffRobin}. If scalar, apply the same value at all vertices. Default: 1.
 \itemtt{nodeData} The vector of nonhomogeneous vertex terms $\phi_n$ in the Poisson problem~\eqref{VC}. If scalar, apply the same value at all vertices. Default: 0.
 \itemtt{plotCoordinateFcn} The handle of a function defining the layout of the edges and vertices for plotting. Associates coordinate arrays \listtt{x1}, \listtt{x2}, and (optionally) \listtt{x3} to each edge and to the vertices. If left unset, then plotting is not possible. It can be set later using the function \listtt{addPlotCoordinates}.
 \end{description}
 The constructor runs several checks on the inputs to ensure they are consistent and meaningful, returning descriptive error messages if these checks fail.

\subsection{Properties of a \qGraph\ object}
Many of the properties of a designated \qGraph\ object are detailed in Sec.~\ref{sec:initializing}, a complete list is given here, filling in some additional details
\begin{description}
\itemtt{qg} The \digraph\ object, consisting of \Edge\ and \listtt{Node} tables, each of which has the additional required fields described in Sec.~\ref{sec:initializing} as well as the optional fields \listtt{x1}, \listtt{x2}, and \listtt{x3} used for plotting.
\itemtt{discretization} A string labeling the discretization type is used to choose between uniform and Chebyshev algorithms.
\itemtt{wideLaplacianMatrix} The Laplacian matrix $\Lint$, with discretized boundary condition rows at the bottom, defined in Eq.~\eqref{Luniform} and illustrated by the two upper matrix blocks in Figs.~\ref{fig:lasso}(b) and~\ref{fig:lassocheb}(b).
\itemtt{interpolationMatrix} The matrix $\Pint$ that interpolates from the extended grid to the interior grid as defined in Eq.~\eqref{Puniform}, as illustrated by the two upper matrix blocks in  Figs.~\ref{fig:lasso}(c) and~\ref{fig:lassocheb}(c).
\itemtt{discreteVCMatrix} The matrix $\MVC$ containing the discretization of the vertex conditions, as defined in Eq.~\eqref{Luniform}, \eqref{Puniform} and illustrated by the two lower matrix blocks in Figs.~\ref{fig:lasso}(b) and~\ref{fig:lassocheb}(b).
\itemtt{nonhomogeneousVCMatrix} The matrix $\MNH$ defined in Eq.~\eqref{MNH} used to define nonhomogeneous terms in the vertex condition to the correct rows.
\itemtt{derivativeMatrix} The square first derivative matrix which does not include boundary conditions. This is used for calculating integrals, including the energy and momentum, which may or may not be conserved based on the vertex conditions.
\end{description}

\subsection{Methods defined for a \qGraph\ object}

\subsubsection{MATLAB \digraph\ methods overloaded for \qGraph\ objects}
MATLAB features many functions for analyzing, querying, and manipulating directed graphs. The command \listtt{indegree(G,1)} returns the incoming degree of the vertex $\vertex_1$ of a graph \ttG. This could be applied to the \listtt{qg} field of a quantum graph $\Phi$ by using the command \listtt{indegree(G.qg,1)}, but it is preferable in object-oriented programming to \emph{overload} this function so that can be applied directly as \listtt{indegree(Phi,1)} Several other low-level directed graph functions have been similarly overloaded: 
\begin{description}
\itemtt{Edges, Nodes, indegree, outdegree, numedges, numnodes, rmnode.}{}
\end{description}

\subsubsection{Other \qGraph\ methods}
The following provide directed graph related functionality not in MATLAB's digraph toolbox:
\begin{small}
\begin{description}
\itemtt{source, target, follows, sharednode, incomingedges, outgoingedges, isleaf.}{}
\end{description}
\end{small}
The following functions query specific properties of quantum graphs, edges, or vertices:
\begin{description}
\itemtt{nx, dx, weight, L, robinCoeff, isUniform, isChebyshev, isDirichlet.}{}
\end{description}
The following are utilities for working with \qGraph\ objects:
\begin{description}
\itemtt{addPlotCoords} Given a user-provided script defining the plotting coordinates \listtt{x1}, \listtt{x2}, and, optionally, \listtt{x3}, runs the script and associates the coordinates to both the edge and vertex tables.
\itemtt{graph2column} and \listtt{column2graph} transfer data back and forth between the edge-vertex representation and a single-column vector. The latter function uses the discretized vertex conditions to interpolate the data at the vertices. 
\itemtt{applyFunctionToEdge} The call \listtt{G.applyfunctionToEdge(fhandle,m)} applies the function represented by the function handle \listtt{fhandle} to the edge $\edge_m$ and stores the result in \listtt{G.Edges.y\{m\}}. If \listtt{fhandle} is a number \listtt{c}, then the output \listtt{G.Edges.y\{m\}} will be a constant-valued vector of the appropriate length.
\itemtt{applyFunctionsToAllEdges} If \listtt{handleArray} is a cell array containing $\abs*{\EE}$ function handles and constants, this function applies \listtt{applyFunctionToEdge} to each function/constant and edge in the quantum graph. If an output argument is specified, then \texttt{graph2column} is used to assign the function to a column vector.
\item{addPotential} Adds a potential to the graph structure using the same syntax as \listtt{applyFunctionsToAllEdges}.
\itemtt{applyGraphicalFunction} This applies a function, input as its function handle, to the plotting coordinates \listtt{x1}, \listtt{x2}, and (optionally) \listtt{x3} defined for each function and edge. This convenience function creates initial guesses for the nonlinear standing wave solvers.
\itemtt{addEdgeField} and \listtt{addNodeField} can be used to assign other fields to the \listtt{Edge} and \listtt{Node} tables.
\end{description}
The following functions perform mathematical operations on \qGraph objects, automatically choosing the appropriate program for the discretization method used:
\begin{description}
\itemtt{integral} Computes the weighted integral $ \int_\Gamma \Psi \dd{x} = \sum_{m=1}^{\abs{\EE}} w_m\int_{\edge_m}\psi_m(x) \dd{x}.$
\itemtt{norm} Uses \listtt{integral} to compute the $L^p$ norm~\eqref{Lpnorm}.
\itemtt{dot} Uses \listtt{integral} to compute the $L^2$ inner product~\eqref{innerProduct}.
\itemtt{energyNLS} Uses \listtt{integral} to compute the NLS energy~\eqref{energyNLS}.
\itemtt{eigs} Computes $n$ eigenvalues closest to zero.
\itemtt{secularDet} Computes the real-valued secular determinant defined briefly in Sec.~\ref{sec:eigenvalueProblem} using the MATLAB Symbolic Mathematics Toolbox. This works for all the boundary conditions discussed in this article but requires the edge weights to satisfy $w_m \equiv 1$.
\itemtt{solvePoisson} Solves the Poisson problem~\eqref{Poisson}.
\end{description}
The following functions are for visualizing \qGraph\ objects:
\begin{description}
\itemtt{plot} The call \listtt{G.plot} plots the data currently stored in the \tty entries of the \listtt{Edges} and \listtt{Nodes} tables, using the coordinates stored in the \listtt{x1}, \listtt{x2} and \listtt{x3} table entries. If \listtt{x3} is not defined, then it plots the function in three dimensions over the skeleton of the graph. If it is defined, then the function is plotted in false color. The call \listtt{G.plot(z)} first calls \listtt{G.column2graph(z)} and then plots.
\itemtt{pcolor} Plots the function in false color on the quantum graph in two dimensions. It is useful for visualizing highly complex graphs, as seen by comparing the two plots of MATLAB's \listtt{peaks} function defined over the edges of a randomly generated Delaunay triangulation, shown in Fig.~\ref{fig:delaunay}. 

\begin{lstlisting}[numbers=left,firstnumber=1,numberstyle=\tiny]{matlab}
G = delaunaySquare('n',8);
f = @(x1,x2)peaks(6*x1-3,6*x2-3);
G.applyGraphicalFunction(f);
G.plot; figure; G.pcolor
\end{lstlisting}
\begin{figure}[htbp] 
   \centering
   \includegraphics[width=0.8\textwidth]{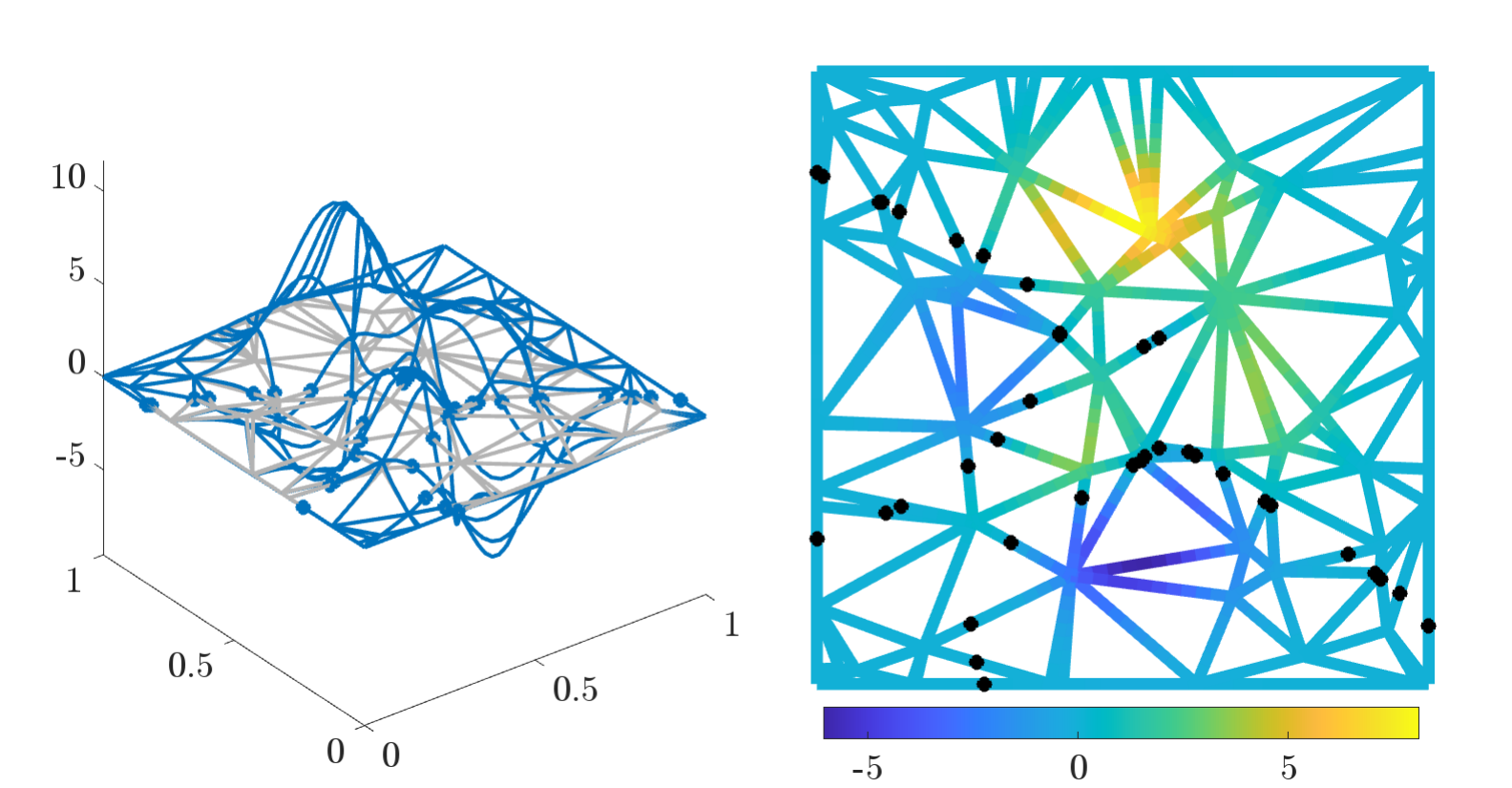} 
   \caption{Visualization of a function defined on a random graph using (left) \listtt{plot} and (right) \listtt{pcolor}, where zeros are indicated with black dots.}
   \label{fig:delaunay}
\end{figure}

\itemtt{spy} The call \listtt{G.spy} uses the MATLAB \listtt{spy} function to plot the nonzero entries in the three matrices \listtt{G.wideLaplacianMatrix}, \listtt{G.interpolationMatrix}, and \\ \listtt{G.nonhomogeneousVCMatrix}.
\itemtt{animatePDESolution} Given a vector of times \listtt{t} and an array \listtt{u} whose columns give the numerical solution to a PDE at those times, animates the solution, taking special care that the viewing axes are fixed throughout the visualization. Automatically uses false color to plot graphs with a three-dimensional layout. To animate a PDE solution using false color on a two-dimensional layout, use \listtt{animatePDESolution2DColor}.
\end{description}
Additional programs not called by the end-user exist, which we do not document.

\subsection{The template library}
The package features a library of graphs, many of which have been studied in the quantum graph literature, which is stored in the folder \listtt{source/templates}. Their use is demonstrated in the live script that is titled \listtt{templateGallery.mlx}. These fall into a few groups. Almost all depend on several user-provided parameters for which default values are provided.  

\subsubsection*{Individual graphs}
\begin{figure}[htbp] 
   \centering
   \includegraphics[width=3in]{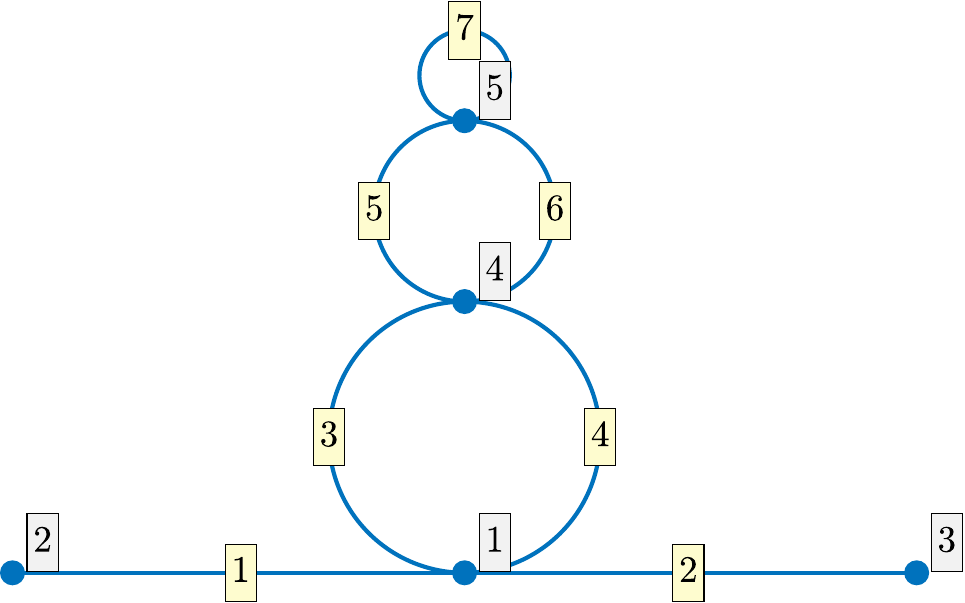} 
   \caption{The default \listtt{bubbleTower} quantum graph, with five vertices and seven edges.}
   \label{fig:bubbletower}
\end{figure}

Several simple graphs are provided in the template library and are called using the command 
\listtt{G=quantumGraphFromTemplate(tag,varargin)}, where \listtt{tag}\ is the name of the template and \listtt{varargin} is used by MATLAB to indicate a variable-length input argument list, and is here used to enter using the same key-value syntax as the \qGraph\ command. The graph produced by running:\\
\begin{footnotesize}
\noindent\listtt{G=quantumGraphFromTemplate('bubbleTower','L',10,'circumferences',[6 4 2]*pi)}
\end{footnotesize}

\noindent is shown in Fig.~\ref{fig:bubbletower}.
The default graph in this family has five vertices and seven edges.  Bubble tower graphs with infinite-length base edges have featured extensively in the quantum graph literature as examples where one can still find a ground state even though a certain graph topology condition is satisfied by this family that would normally preclude the existence of a ground state, see~\cite{adami2016ground,adami2017negative,adami2017nonlinear}.  The underlying symmetry of the construction here is crucial to the analysis.

The \listtt{quantumGraphFromTemplate} function calls two separate functions
\begin{itemize}
\item \textbf{A template function}, here \listtt{bubbleTower.m}, that builds the quantum graph, setting the lengths of the two straight line segments to 10 and the circumferences of the three bubbles to $[6\pi, 4\pi, 2\pi]$, setting the discretization, and building the necessary matrices.
\item \textbf{A plot coordinates function}, here \listtt{bubbleTowerPlotCoords.m}, that places the vertices at locations consistent with the above-defined lengths. In this example, two edges are laid out as line segments, created using the command \listtt{straightEdge}, four edges are laid out as semicircular edges by using the command \listtt{semicircularEdge}, and there is one circular edge, created using the command \listtt{circularEdge}. A fourth function \listtt{circularArcEdge} can connect two nodes by a circular arc subtending a central angle \listtt{theta}.
\end{itemize}

\subsubsection{Two-dimensional lattices}
The following templates exist to create two-dimensional lattices. All have default values and can be customized to change the number of cells per side.  These programs are called directly and set plotting coordinates without calling \listtt{quantumGraphFromTemplate}.
\begin{description}
\itemtt{rectangularArray} creates a rectangular array. By default, the sides have unit length but can be customized.
\itemtt{triangularArray} creates a triangular array. The unit cell is an equilateral triangle by default, but the period vectors can be customized.
\itemtt{hexagonalArray} creates a hexagonal array, forming a parallelogram, the default shown in Fig.
\itemtt{hexGrid} creates a rectangular array of hexagons.
\itemtt{hexGridPeriodic} identifies the left edge with the right and the top edge with the bottom to create a periodic array.
\itemtt{hexOfHexes} A hexagonal array of hexagons.
\itemtt{triangularArray} A triangular array.
\end{description}

\subsubsection*{Three-dimensional geometric templates}
The program \listtt{solidTemplate} constructs quantum graphs whose vertices and edges are the vertices and edges of geometric solids, including the five Platonic solids (tetrahedron, cube, octahedron, dodecahedron, and icosahedron), as well as the cuboctahedron, which has 24 edges and 12 vertices, and the buckyball (or truncated icosahedron) which has 90 edges and 60 vertices. This is called directly and sets up the plot coordinates. Sec.~\ref{sec:functions-and-plotting} gives an example of constructing a tetrahedron.

\subsection{Continuation and bifurcation routines} 
\label{sec:appendixContinuation}
The live script that is titled \listtt{continuationInstructions.mlx} in the \listtt{documentation} directory uses all the following subroutines in the given order after constructing a \qGraph\ object named \listtt{Phi}. We refer to line numbers in this live script to describe the steps taken to compute the bifurcation diagrams. To run the continuation software, the user must use a template from the \listtt{source/templates} or create one themselves, including a properly named function to create the plotting coordinates. We will assume that the template's name is stored in a variable named \listtt{tag}. In the example \listtt{tag='dumbbell'}. As explained below, the results of the computation will be stored in the directory \listtt{dataDir='data/dumbbell/001'} with the trailing number incremented each time a bifurcation diagram is created. Each computed branch of solutions is stored in its own subdirectory, with consecutively labeled names, beginning \listtt{branch001}, etc. Most of the programs given below add a line to a log file named \listtt{logfile.txt} that resides in the data directory.

\begin{description}
\itemtt{makeContinuationDirectory} After initializing the discretized quantum graph on which families of solutions are to be computed, create a sequentially named directory to hold the data; see line \listtt{5}. Saves a file \listtt{template.mat} containing the \qg\ object.
\itemtt{saveEigenfunctions} Calculate some eigenvalues and eigenfunctions of the Laplacian matrix and save them to the data directory with names \listtt{lambda.001} and \listtt{eigenfunction.001}. 
\itemtt{saveNLSFunctionsGraph} Saves a file named \listtt{fcns.mat} to the data directory. This file contains one variable: a structure \ttx\ whose fields contain a function handle to the discretized form of Equation~\eqref{stationary_NLS}, as well as several derivatives of this function and the antiderivative of the nonlinearity, used in computing the energy.
\itemtt{continuerSet} This function sets several parameters the continuation algorithms use. It assigns them to a structure, usually named \listtt{options}, which is then passed to the various \listtt{continueFrom} programs described below. It takes as input a sequence of name-value pairs, imitating the programs \listtt{odeset} and \listtt{optimset} used in MATLAB's ODE and optimization routines. The parameters it sets are:
\begin{description}
\itemtt{maxTheta} The maximum angle, in degrees, between two consecutive segments on a branch of solutions. Default: $4^\circ$.
\itemtt{minNormDelta} The minimum step length below which the continuation solver does not attempt to refine the branch further. Default: $10^{-3}$.
\itemtt{beta} The weight in the inner product defined by  \[
\ip{\Phi_1(x) e^{i \Lambda_1 t}}{\Phi_2(x) e^{i \Lambda_2 t}} = \ip{\Phi_1}{\Phi_2} + \beta \ip{\Lambda_1}{\Lambda_2},
\] 
used in defining angles and distances in the above two variables. Default: $0.1$.
\itemtt{NThresh} Threshold for the power $N$, i.e., the squared $L^2$-norm, so the continuation routine terminates when this value is crossed. Default: 4. 
\itemtt{LambdaThresh} Threshold for the frequency $\Lambda$. The continuation routine terminates when this value is crossed. Default: -1. 
\itemtt{maxPoints} The maximum number of points to compute on a given branch. Default: 999.
\itemtt{saveFlag} A boolean variable. If true, then data is saved to files. Default: \listtt{true}.
\itemtt{plotFlag} A boolean variable. If true, then data is plotted to screen. Default: \listtt{true}.
\itemtt{verboseFlag} A boolean variable. If true, then some information is printed on the MATLAB Desktop. Default: \listtt{true}.
\end{description}
\item Four continuation programs that are initiated from different starting points.
\begin{description}
\itemtt{continueFromEig} Compute a branch of stationary solutions that bifurcates from $\Psi=0$ with a frequency given by an eigenvalue, in the direction of an eigenfunction, using the data saved by the above command \listtt{saveEigenfunctions}; cf.\ lines \listtt{13-15} of the live script.
\itemtt{continueFromBranchPoint} Compute a branch of stationary solutions that bifurcates from a branch point. While computing a curve of solutions, the continuation routines monitor for branching bifurcations (pitchfork and transcritical, which are mathematically equivalent in the pseudo-arclength formulation). When it detects a bifurcation between two computed solutions, it computes the exact frequency at which the bifurcation occurs and the solution at the bifurcation point.
\itemtt{continueFromSaved} Continue from a previously-computed solution to the stationary computed using \listtt{saveHighFrequencyStandingWave} (called here), which computes and saves a solution with an initial guess built from $\sech$-like functions defined on the edges, \\ \listtt{saveHighFrequencyStandingWaveGraphical}, \\
which computes a solution based on an initial guess that places a "bump" somewhere on the graph defined by its plotting coordinates or a user-written function. 

On line \listtt{35} of the example, a solution with positive $\sech$ pulses of edges 1 and 2 of the dumbbell is saved to files: \\
\listtt{savedFunction.001} and \listtt{savedFunction.001} \\
in the folder \listtt{data/dumbbell/001}. A branch continuing from this solution is computed at line \listtt{38}.
\itemtt{continueFromEnd} Extends a previously-computed branch.
\end{description}
\itemtt{bifurcationDiagram} Draws a bifurcation diagram from the data in a given directory and its subdirectories. By default, it plots the frequency on the $x$-axis and the squared $L^2$-norm on the $y$-axis, but these defaults can be overwritten.
\itemtt{rmBranch} Removes the subdirectory containing a given branch from the bifurcation diagram directory.
\itemtt{plotSolution} Plots a single solution from a given diagram and branch.
\itemtt{animateBranch} Animates how the individual solutions change as a branch of the bifurcation diagram is traversed. 
\itemtt{addComment} Adds a string to the log file \listtt{logfile.txt} in the given directory.
\end{description}

We examine the files contained in the directory \listtt{branch001}, which was created online \listtt{13} of the live script by \listtt{continueFromEig}. 
\begin{description}
\itemtt{PhiColumn.xxx} Where \listtt{xxx} is a three-digit integer $n$. The $n$th solution on the branch.
\itemtt{NVec}, \listtt{LambdaVec}, and \listtt{energyVec} Column vectors containing the squared $L^2$-norm, the frequency, and the energy, which are the three variables that can be plotted using the \listtt{bifurcationDiagram} program. The $n$th entry in each vector corresponds to the $n$th solution in the previous bullet point.
\itemtt{k} The number of \listtt{PhiColumn} files and the length of the vectors of integrals.
\itemtt{initialization} A one-word text file denoting which of the four continuation programs \listtt{coninueFromXXX} was used to initialize the branch, in this case \listtt{Eigenfucntion}.
\itemtt{eignumber} The number of the eigenfunction from which the solution was continued.
\itemtt{options.mat} The options structure set by the \listtt{continuerSet} program.
\itemtt{bifTypeVec} A column vector of integers, with the value 0 if solution $n$ is a regular point on the branch, the value 1 at branching bifurcations, and the value -1 at folds.
\itemtt{phiPerturbationXXX.mat} and \listtt{LambdaPerturbationXXX.mat} Here \listtt{xxx} is a three-digit number at which a branching bifurcation has been detected, and the files contain the directions in which the new branch points from the bifurcation location, used by the function \listtt{continueFromBranchPoint}.
\end{description}

\subsection{Other folders}
\begin{description}

\itemtt{data} An empty folder where the continuation routines store the data they produce.
\itemtt{documentation} Contains live scripts demonstrating the main features entitled  \\
\listtt{quantumGraphRoutines.mlx}, \listtt{continuationInstructions.mlx}, and \\ \listtt{continuationInstructionsChebyshev.mlx}.
\itemtt{source/chebyshev} Contains many programs used to construct the Chebyshev discretization.
\itemtt{source/examples} Contains example programs sorted into three further subfolders:
\begin{description}
\itemtt{source/examples/chebyshev} Contains examples involving  the Chebyshev discretization, all of which are minor modifications of examples from the \listtt{stationary} folder.
\itemtt{source/examples/evolution} Examples illustrating the solution to time-dependent problems.
\itemtt{source/examples/stationary} Examples of time-independent problems: eigenproblems, Poisson problems, and continuation problems.
\end{description}
\itemtt{source/startup\_shutdown} Contains programs that are run upon starting up and shutting down QGLAB.
\itemtt{source/user} An empty folder intended to give end-users a place to store code they write without mixing it with package code.
\itemtt{source/utilities} Some utilities used for file management and formatting plots.
\itemtt{tmp} A temporary folder created at startup and removed at shutdown, where the user's plotting preferences are stored to be automatically restored upon shutting down \qGraph.
\end{description}

\section*{Acknowledgments}
The authors thank Greg Berkolaiko, Dmitry Pelinovsky, David Shirokoff, and Nick Trefethen for their helpful conversations.  RG credits a semester-long funded visit to the IMA at the University of Minnesota in 2017 for giving him the time and inspiration to begin working on quantum graphs.

\end{document}